\documentclass[a4paper,11pt]{article}
\usepackage{amssymb}
\usepackage{amsthm, amsmath, amscd}
\newtheorem{lemma}{\bf Lemma}[section]%
\newtheorem{theorem}[lemma]{\bf Theorem}%
\newtheorem{proposition}[lemma]{\bf Proposition}% 
\newtheorem{corollary}[lemma]{\bf Corollary}%   

\newtheorem{remark}[lemma]{\sc Remark}% 

\newtheorem{definition}[lemma]{\bf Definition}{}%  
\newcommand{\be}{\begin{equation}}
\newcommand{\ee}{\end{equation}}
\newcommand{\ol}{\overline}

\newcommand{\wh}{\widehat}

\begin{document}
%\begin{frontmatter}


 \title{Edge colouring models for the Tutte polynomial and related
 graph invariants}
% \thanks[label1]{Partially supported by EPSRC}
 \author{A.J. Goodall} 

% \ead{a.goodall@bristol.ac.uk}
% \address{University of Bristol, University Walk, Bristol BS8 1TW, United Kingdom.}
%\date
\maketitle

\begin{abstract}
For integer $q\geq 2$, we derive edge $q$-colouring models for (i) the
Tutte polynomial of a graph $G$ on
the hyperbola $H_q$, (ii) the
symmetric weight enumerator for group-valued $q$-flows of $G$, and
(iii) a more general vertex colouring model partition function that includes these
polynomials and the principal
specialization order $q$ of Stanley's symmetric monochrome polynomial. We
describe the general relationship between vertex and edge colouring
models, deriving a result of Szegedy and generalizing a
theorem of Loebl along the way.
In the second half of the paper we exhibit a family of non-symmetric edge
$q$-colouring models defined on $k$-regular graphs, whose partition
functions for $q\geq k$ each evaluate the number of proper edge
$k$-colourings of $G$ when $G$ is Pfaffian.    
\end{abstract}

%\begin{keyword}
%graph \sep edge colouring model \sep vertex colouring model \sep chromatic polynomial \sep Tutte polynomial \sep symmetric weight
%enumerator \sep symmetric chromatic function \sep  discrete Fourier transform \sep
%proper edge colouring

%\PACS  Coloring of graphs and hypergraphs \sep Factorization, matching, covering and packing
%\sep graph \sep partition
%function \sep edge colouring model \sep vertex colouring model \sep Fourier transform \sep proper edge colouring
% \PACS 05C15 \sep 05C70
%\end{keyword}

%\end{frontmatter}
\tableofcontents
%\listoffigures

\section{Introduction}

\subsection{Background}\label{sec: background}

Partition functions of edge colouring and vertex colouring models are
graph invariants with a definition motivated by statistical physics,
where the vertices or edges of a graph $G$ are ``sites'' and colours are ``states''. The states
interact with each other along the edges in vertex colouring models
and at vertices in edge colouring models. A configuration (assignment
of states to all sites) is given a weight dependent on the
interactions between states and possibly also on the states
themselves. The partition function of the model is obtained by summing
the weights over all possible configurations. 
For example, the partition function of
the $q$-state Potts model \cite{Potts52} is defined by a vertex
                                colouring model, where the interaction
                                along an edge depends on whether
                                the endpoint states are the same or
                                different. 

Freedman, Lov\'{a}sz and Schrijver \cite{FLS07} have shown that a
graph parameter can be realized as the partition function of a
real-valued vertex colouring model if and only if it is vertex reflection
positive and has exponential rank-connectivity. Szegedy \cite{S05}
proved the parallel statement, that graph parameters arise as
partition functions of real-valued edge colouring models if and only
if they are edge reflection positive and multiplicative over
disjoint graphs.
(To be reflection positive means, approximately, to take non-negative
values on graphs with mirror symmetry.)  Vertex and edge
colouring models thus play an important r\^{o}le in the algebra of quantum graphs, for
which see the papers already cited and for example \cite{BCLSV}, \cite{LS}, \cite{LCCC07}, \cite{L06}. 

Schrijver  \cite{LS07}, in the more general context of hypergraphs, extends the
scope of vertex colouring models to directed graphs (where the
order in which colours appear on an edge now matters) and edge
colouring models to graphs where the edges incident with a vertex come
in a given order. Examples of partition functions of vertex colouring models on
directed graphs include the number of oriented
colourings, or the number of weak oriented colourings (monochromatic
edges allowed) that feature in Stanley's \cite{St73} evaluation of the chromatic polynomial at negative integers.

An embedding of a graph on an orientable surface is defined by its
 vertex rotations (a clockwise order of incident edges as seen from a
 fixed side of the surface).
For embedded graphs it would be natural to consider edge colouring models that depend on the order up to cyclic
 permutation in which colours appear at a vertex. Given a $4$-regular
 plane graph with a chequerboard colouring of its faces, the Viergruppe (elementary $2$-group order $4$) fixes the property of vertex transitions being black, white or
crossing. Edge colouring models with this symmetry might feature
in partition functions for evaluations of a transition polynomial. (For
transition polynomials, of which the Penrose polynomial
is an instance, see for example \cite{MA97}, \cite{MA01}, \cite{FJ89}, \cite{FJ90}, \cite{Sar05}.)

Edge colouring models that depend only on the order up to even
permutation in which colours appear at a vertex
play a r\^{o}le in proper edge $k$-colouring
$k$-regular Pfaffian graphs \cite{EG05}, \cite{NT07}. The overall parity of
the permutations determined by the distinct colours appearing at each vertex turns
out to be fixed for such graphs.  This property allows us to derive
some non-symmetric edge $q$-colouring models with a partition function that
gives the number of proper edge $k$-colourings of the graph.

\subsection{Outline}

The first purpose of this article is to give a general account of
vertex and edge colouring models and their interrelation, 
and in particular to provide examples of such models for some
important graph polynomials, supplementing those given for example in \cite{BCLSV}, \cite{LCCC07}, \cite{S05}. 
 
The second purpose of this article is to present some examples of edge colouring models
that depend on the order in which colours appear around a vertex, while still giving a partition function which is an interesting graph
parameter. 

We begin in Section \ref{sec: prelim} by introducing relevant
notation and concepts. 

In Section \ref{sec: vertex models edge} edge colouring models are
derived for a family of
partition functions of vertex colouring models, including two
branches descending from the chromatic polynomial. This family includes the Tutte polynomial\footnote{See for example
\cite{DW93}, \cite{DW99}, \cite{WM00}, \cite{Sokal} for an account of the significance of
this polynomial in combinatorics.} in its specialization to the
$q$-state Potts model, treated in Section \ref{Tutte} along with its
generalization to the symmetric weight enumerator of
abelian-group-valued $q$-flows (Theorem
\ref{thm: cwe vertex edge}). Szegedy's result \cite{S05} that any symmetric real-valued vertex
colouring model has an edge colouring model is included in
Section \ref{Vertex to Edge}.  The family described in Section \ref{hierarchy} includes the principle
specializations of finite order of Stanley's symmetric monochrome
polynomial \cite{St98} (the latter a generalization of his symmetric chromatic function \cite{St95}). Theorem
\ref{thm: van der W} is a generalization of a theorem of Loebl
\cite{Loebl06}, itself a generalization of van der Waerden's eulerian
subgraph expansion of the Ising model \cite{vdW41}. We finish in Theorem \ref{thm:
  general chrom edge} with an edge colouring model that unifies the aforementioned polynomials.

In Section \ref{Cyclic} we derive the already advertised example of a family of
non-symmetric edge colouring models. The partition function of the edge $q$-colouring model given in Theorem
\ref{thm: PLGk even odd} evaluates a particular coefficient of
Petersen's graph monomial \cite{Pet1891} of the line graph of $G$, considered more
recently by Alon
and Tarsi \cite{AT92}, Matiyasevich \cite{M04} and Ellingham and
Goddyn \cite{EG05}. For the restricted class of
$k$-regular graphs that admit Pfaffian labellings \cite{NT07}, \cite{RT07} (a class that includes
planar graphs and the Petersen graph), this coefficient is up to sign
equal to the number of proper edge $k$-colourings of $G$.
The two main theorems are to be found in Section \ref{sec: last} after the
preparatory work of Sections \ref{sec: 2-factorizations} and \ref{sec:
  fourier parity}.

\subsection{Cubic graphs}

The main theorems of this paper are most easily described at this stage by
giving their statements for a cubic ($3$-regular) graph $G$. 
These special cases give
the flavour of the more substantial results that are proved in the body
of the paper without the need for much preparation or notation. 
For the purpose of illustration we shall
focus on $F(G;4)$, the number of nowhere-zero
$\mathbb{F}_4$-flows (or $\mathbb{Z}_4$-flows) of a cubic graph $G$.

Let $Q$ be an abelian group of order $q$. For a graph $G$, the number
of nowhere-zero $Q$-flows is independent of the structure of $Q$,
equal to $F(G;q)$, the flow polynomial evaluated at $q$.
It is well known (see e.g. \cite{DW93}) that, for any graph $G$,
$$F(G;q)=(-1)^{|E|}q^{-|V|}\sum_{\mbox{\rm \tiny vertex $q$-colourings
    of $G$}}(1-q)^{\#\mbox{\rm \tiny monochrome edges}},$$
where an edge is monochrome if its endpoint vertices have the same
colour. This is an example of the partition function of a vertex $q$-colouring model.

A consequence of Corollary \ref{cor: Tutte vertex edge} (by setting $s=0$) is that $F(G;q)$ is also the partition
function of an edge colouring model, which for cubic $G$ has a
particularly concise expression:\medskip

\begin{proposition} \label{prop: q-flows cubic edge}For a $3$-regular graph $G=(V,E)$ and $q\geq 2$,
  $$F(G;q)=q^{-|E|}2^{|V|}\sum_{\mbox{\rm \tiny edge $q$-colourings of $G$}}(1-q)^{\#\,\mbox{\rm \tiny  monochrome
    vertices}}(1-q/2)^{|V|-\#\mbox{\rm \tiny rainbow vertices}},$$
where a vertex is monochrome (rainbow) in an edge $q$-colouring if
the colours on its
incident edges are the same (different).  
\end{proposition}

Theorem \ref{thm: cwe vertex
  edge} gives an edge $q$-colouring model for a generalization of
  $F(G;q)$, namely the symmetric weight enumerator of the set of $Q$-flows of
  $G$. (The Tutte polynomial on the hyperbola $H_q$ is the Hamming
  weight enumerator of the set of $Q$-flows. An edge $q$-colouring
  model for this is given in  Corollary \ref{cor: Tutte vertex edge}.)  

Vertex colouring models are uniquely determined \cite{S05}, in the sense that
different vertex colouring models give partition functions that differ
for some graph $G$. Edge colouring models on the other hand may be
rotated by an orthogonal transformation and yet preserve the partition
function (Szegedy \cite{S05}, and Theorem \ref{thm: szegedy}
below). A reverse phenomenon is seen to occur here, however, for when $G$
is cubic $F(G;4)$ has an infinite
number of vertex $4$-colouring models all of which, on applying
Theorem \ref{thm: cwe vertex edge}, come from the same
edge colouring model ($q=4$ in Proposition \ref{prop: q-flows cubic
  edge}).\medskip

\begin{proposition}\label{prop: Fg4 vertex}
For a $3$-regular graph $G=(V,E)$ and indeterminates $s,t$,
$$(st)^{|E|/3}F(G;4)=4^{-|V|}\sum_{\mbox{\rm \tiny vertex
    $\mathbb{F}_4$-colourings}}(1+s+t)^{\#
  0}(1-s-t)^{\# 1}(-1-s+t)^{\# \omega}(-1+s-t)^{\# \ol{\omega}},$$
where the exponent $\# a$ for $a\in\mathbb{F}_4=\{0,1,\omega,\ol{\omega}\}$ counts the number of
edges whose endpoint colours differ by $a$. 
\end{proposition}
[This proposition is not proved below, but is a simple application of
MacWilliams duality identity (given as equation \eqref{eqn: Poisson}
later) on observing that all nowhere-zero $\mathbb{F}_4$-flows of $G$
have an equal number of each non-zero element.]

The vertex colouring model in Proposition \ref{prop: Fg4 vertex} depends on the structure of $\mathbb{F}_4$
in order to define the weight it gives to an edge in a vertex
$4$-colouring. It loses this dependence on setting $s=t=u/2$, when it can
be written as a vertex $\mathbb{Z}_4$-colouring model: 
$$u^{|V|}F(G;4)=(-1)^{|E|}2^{|V|}\sum_{\mbox{\rm \tiny vertex $\mathbb{Z}_4$-colourings}}(-1-u)^{\# 0}(-1+u)^{\#\, 2},$$
where $\# 0$ is the number of monochrome edges and $\# 2$ denotes the
number of edges with colour difference $2$.

So far all our vertex and edge colouring models have been
symmetric.  
They do not depend on the order in which colours appear on the
endpoints of an edge (for vertex colouring models) or the order of
colours appearing on the edges incident with a vertex (for edge
colouring models). This is common to all the vertex colouring models
for which edge colouring models are derived in Section \ref{sec: vertex models edge}. 

In Section \ref{Cyclic} we consider edge
colouring models that depend on the order of colours up to
even permutation. For cubic graphs this coincides with order up to
cyclic permutation, and this has a natural interpretation in terms of
clockwise and anticlockwise rotations. 
 
Suppose a cubic graph $G$ is embedded in an orientable surface. Proper
edge $3$-colourings of $G$ are the same as nowhere-zero
$\mathbb{F}_4$-flows of $G$, and from this comes Tait's
\cite{Tait1880}, \cite{Tait80} equivalent
statement of the Four Colour Theorem, that every planar cubic graph has
a proper edge $3$-colouring. By Vizing's theorem \cite{Vizing64} every
$k$-regular graph has a proper edge $(k+1)$-colouring. Theorem
\ref{thm: k+1} relates proper edge $k$-colourings to proper edge
$(k+1)$-colourings of $k$-regular graphs, its special case for plane
cubic graphs being the proposition that follows here. Theorem
\ref{thm: PLGk even odd} more generally gives for any $q\geq k$ an edge $q$-colouring model for the number of proper edge
$k$-colourings of a $k$-regular graph that admits a Pfaffian labelling \cite{NT07}, \cite{RT07}.

An orientable
embedding of $G$ is described by its vertex rotations, giving a
clockwise order of edges around a vertex. 
Let $0<1<2<3$ be ordered up to cyclic permutation, i.e. as the
cycle $(0\; 1\; 2\; 3)$. Given a proper edge
$4$-colouring of $G$ with colours $\{0,1,2,3\}$, the three colours
 that appear at a vertex come either in a clockwise sense, i.e. consistently with the cyclic order $(0\; 1\; 2\;
3)$, or in an anticlockwise sense, i.e. in the reverse order  $(3\;
2\; 1\; 0)$. Colours can appear in a clockwise order
(e.g. $0,2,3$), or in an anticlockwise order (e.g. $0,3,2$).\medskip
 
\begin{proposition} \label{prop: cubic even odd 4} Let $G=(V,E)$ be a plane cubic graph. Then
$$(-4)^{|E|/3}F(G;4)$$
$$=\#\{\mbox{\rm \small even proper edge
  $4$-colorings of $G$}\}-\#\{\mbox{\rm \small odd proper edge
  $4$-colorings of $G$}\},$$
where a proper edge $4$-colouring of $G$ is even (odd) if there are an
even (odd) number of vertices at which the colours appear in
an anticlockwise order. 
\end{proposition}

\section{Preliminaries}\label{sec: prelim}

For more on elementary Fourier analysis
see for example \cite{Terras99}. 
The assumed graph theory is standard.
For boundaries and coboundaries, and also for a related perspective on vertex
colouring models, see \cite{B77}.  

\subsection{Abelian groups and the Fourier transform}\label{sec:
  groups fourier}
Let $Q$ be a finite additive abelian group order $q$, which we assume also has
the structure of a commutative ring with unity. 

The set $\mathbb{C}^Q$ of all functions $f:Q\rightarrow\mathbb{C}$ is
an inner product space with Hermitian inner product
$$\langle f,g\rangle=\sum_{a\in Q}f(a)\ol{g(a)},$$
where the bar denotes complex conjugation. The subspace $\mathbb{R}^Q$
has Euclidean inner product 
$$(f,g)=\sum_{a\in Q}f(a)g(a),$$
and we may also use the notation $(\,,\,)$ in the larger space
$\mathbb{C}^Q$. %, noting that $(f,g)=\langle f,g\rangle$ if $g\in\mathbb{R}^Q$.

The pointwise product of $f$ and $g$ is defined by $f\cdot
g(a)=f(a)g(a)$ and their convolution by $f\ast g(a)=\sum_{b\in
  Q}f(a-b)g(b)$.

Elements of $\mathbb{C}^Q$ will be regarded interchangeably as
functions and as column vectors indexed by elements of $Q$.
The indicator function of $P\subseteq Q$ is denoted by $1_P$.

A character of $Q$ is a homomorphism $\chi:Q\rightarrow\mathbb{C}^\times$
from $Q$ to the multiplicative of $\mathbb{C}$. The set of
characters $\wh{Q}$ forms a group isomorphic to $Q$. A character $\chi$ is a {\em generating
  character} for $Q$ if the isomorphism $Q\rightarrow\wh{Q}$ can be
realized by the map $a\mapsto\chi_a$ defined by $\chi_a(b)=\chi(ab)$. Given
a generating character $\chi$, the matrix $$F:=q^{-1/2}(\chi(ab))_{a,b\in
  Q}$$ is the Fourier transform on the space $\mathbb{C}^Q$. The
following are three key properties of the Fourier transform:

\begin{itemize}
\item $F$ is a unitary matrix, i.e. $\ol{F}^TF=I$, so that 
$\langle Ff,Fg\rangle=\langle f,g\rangle$. 

\item Pointwise products are transformed by $F$ into convolutions,
  $F(f\cdot g)~=q^{-1/2}Ff\ast~ Fg.$

\item If $P$ is a subgroup of $Q$ and $1_P$ is the indicator function
  of $P$ then $F1_P=|P|1_{P^\sharp}$, where $P^{\sharp}=\{a\in
  Q:\forall_{b\in P}\;\chi(ab)=1\}$. 
\end{itemize}
The $d$-fold Cartesian
product $Q^d$ is itself an abelian group and a module over
$Q$. Multiplication on $Q^d$ is componentwise and the dot product of elements $a=(a_1,\ldots, a_d)$
and $b=(b_1,\ldots, b_d)$ is defined by $a\cdot b=a_1b_1+\cdots + a_db_d$. If
$\chi$ is a generating character for $Q$ then $\chi^{\otimes d}$,
defined for $(a_1,\ldots, a_d)\in Q^d$ by $\chi^{\otimes d}(a_1,\ldots, a_d)=\chi(a_1)\cdots
\chi(a_d)$, is a
generating character for $Q^d$ and $\chi^{\otimes d}(ab)=\chi(a\cdot b)$. The matrix $F^{\otimes
  d}$ is the
matrix for the Fourier transform on the space $\mathbb{C}^{Q^d}$,
which, as above, satisfies $\langle F^{\otimes d}f,F^{\otimes d}
g\rangle~=\langle f,~g\rangle$ and $F^{\otimes d}(f\cdot g)=q^{-d/2}Ff\ast Fg$ for functions $f,g\in\mathbb{C}^{Q^d}$. 
Unitary transformations are isometries of the Hermitian inner product
space $\mathbb{C}^{Q^d}$, and orthogonal transformations isometries of the Euclidean
inner product space $\mathbb{R}^{Q^d}$. If $U$ is orthogonal, i.e. $UU^T=I$, then $(U^{\otimes
  d}f,U^{\otimes d}g)=(f,g)$ for $f,g\in\mathbb{C}^{Q^d}$.

 For $\mathcal{C}\subseteq Q^d$ the orthogonal submodule to $\mathcal{C}$ is defined by
$$\mathcal{C}^\perp=\{a\in Q^d:\forall_{c\in\mathcal{C}}\;a\cdot c=0\},$$
and,  given that $Q$ has a generating
  character $\chi$, we have
  $\mathcal{C}^{\perp}~=\mathcal{C}^{\sharp}~=\{a\in
  Q^d~:\forall_{c\in \mathcal{C}}~\;\chi(a\cdot c)=1\}$, so that 
$$F^{\otimes d}1_{\mathcal{C}}=q^{-d/2}|\mathcal{C}|1_{\mathcal{C^\perp}}.$$

Define the set of all sequences on $Q$ by $$Q^*=\bigcup_{d\in\mathbb{N}}Q^d,$$ 
and the following subsets:
$$\mbox{\sc Monochrome}=\bigcup_{d\in\mathbb{N}}\{(a_1,a_2,\ldots,a_d)\in
Q^d:a_1=a_2=\cdots=a_d\},$$
$$\mbox{\sc Zero-sum}=\bigcup_{d\in\mathbb{N}}\{(a_1,a_2,\ldots,a_d)\in
Q^d:a_1+a_2+\cdots+a_d=0\}.$$
Note that, for each $d\in\mathbb{N}$,
$\mbox{\sc Zero-sum}\cap Q^d=(\mbox{\sc Monochrome}\cap Q^d)^\perp$,
so that
 $$F^{\otimes d}1_{\mbox{\tiny \sc Monochrome}\cap
   Q^d}=q^{1-d/2}1_{\mbox{\tiny \sc Zero-sum}\cap Q^d}.$$

Finally, it will be useful to have the following notation.
Suppose that $U$ is a linear transformation of
$\mathbb{C}^{Q}$, so that the $d$-fold tensor product $U^{\otimes d}$ is a
linear transformation of $\mathbb{C}^{Q^d}$ for each $d\in\mathbb{N}$. 
Then the map $f\mapsto f^U$ is defined to be the unique algebra homomorphism
$\mathbb{C}^{Q^*}\rightarrow\mathbb{C}^{Q^*}$ satisfying $f^U=Uf$ for
$f\in\mathbb{C}^Q$. In other words, $f^U$ is defined
for all $d\in\mathbb{N}$ and $z\in Q^d$ by $f^U(z)=U^{\otimes
  d}f(z)$. In this way, for example, given a function
$f\in\mathbb{C}^{Q^*}$ the function $f^F$ is the Fourier transform of
$f$ taken in the appropriate space according to the argument of $f$.

\subsection{Graphs and half-edges}
Let $G=(V,E)$ be a graph with set of vertices $V$ and set of edges
$E$. Edges $e\in E$ are subsets of $V$ of size $2$ (a multiset of size
$2$ on $1$ vertex if $e$ is a loop).  Given $e$ we shall write $v\in e$ when $v$ is an endpoint of
$e$, and, given $v$, we shall write $e\ni v$ when $e$ is incident with $v$. The number of connected components of $G$ is denoted by $k(G)$ and its rank
by $r(E)=|V|-k(G)$. For $A\subseteq E$ the rank of the induced subgraph
$(V,A)$ is denoted by $r(A)$. 
 
The set of {\em half-edges} of $G$ is defined by $H=\{(v,e):v\in
e\}$.
There are two natural ways to partition the set of half-edges. The
first is according to incidence with vertices $v\in V$, with blocks
$H(v):=\{(u,e)\in H:u=v\}$
of size $|H(v)|=d(v)$ the degree of $v$.
The second is according to incidence with edges $e\in E$, and here a
block is a set $H(e):=\{(v,f)\in H:f=e\}$ of size $2$. If $e=\{u,v\}$ then $H(e)=\{(u,e),(v,e)\}$ while if
$e=\{v\}$ is a loop then $H(e)$ contains two
copies of the half-edge $(v,e)$.

Since $\{H(e):e\in E\}$ is a partition of $H$ we have the isomorphism
$$\mathbb{C}^{Q^H}\cong \bigotimes_{e\in E}\mathbb{C}^{Q^{H(e)}},$$
where $\mathbb{C}^{Q^{H(e)}}\cong\mathbb{C}^{Q^2}$. In other words, if
$g\in\mathbb{C}^{Q^H}$ then we can write $g=\otimes_{e\in E}g_e$
for functions $g_e\in\mathbb{C}^{Q^{H(e)}}$, i.e., for $z=(z_h:h\in H)\in Q^H$, 
$$g(z)=\prod_{e\in E}g_e(z_h:h\in H(e)).$$
Likewise, since $\{H(v):v\in V\}$ is a  partition of $H$ we also have the isomorphism
$$\mathbb{C}^{Q^H}\cong \bigotimes_{v\in V}\mathbb{C}^{Q^{H(v)}},$$
where $\mathbb{C}^{Q^{H(v)}}\cong\mathbb{C}^{Q^d}$ when $v$ is a
vertex of degree
$d$.

In the other direction, a function $f:Q^*\rightarrow\mathbb{C}$ defines a
function $f^{\otimes V}:Q^H\rightarrow\mathbb{C}$ given by
$$f^{\otimes V}(z)=\prod_{v\in V}f(z_h:h\in H(v)).$$
Similarly, a function $g:Q^2\rightarrow\mathbb{C}$ extends to a function
$g^{\otimes E}:Q^H\rightarrow\mathbb{C}$ defined for $z\in Q^H$ by
$$g^{\otimes E}(z)=\prod_{e\in E}g(z_h:h\in H(e)).$$

\subsection{The boundary and coboundary}

Let $\sigma$ be an orientation of $G=(V,E)$, defined by 
$$\sigma_{v,e}=\begin{cases} +1 & \mbox{$e$ is directed into
    $v$},\\
-1 & \mbox{$e$ is directed out of $v$},\\
0 & \mbox{$v\not\in e$}.\end{cases}$$
%The function $\sigma:V\times E\rightarrow\{0,\pm 1\},\; (v,e)\mapsto \sigma_{v,e}$ may be regarded as a signed indicator
%vector of the set of half-edges $H\subseteq V\times E$, with a
%half-edge $(v,e)$ negative when it is the tail of a directed edge
%and positive when it is the head of a directed edge.
 
The {\em  boundary} operator $\partial:Q^E\rightarrow Q^V$ is defined by
$$(\partial y)_v=\sum_{e\in E}\sigma_{v,e}y_e.$$
The submodule $\ker(\partial)$ is the set of {\em $Q$-flows} of
$G$. 
The {\em  coboundary} operator $\delta:Q^V\rightarrow Q^E$ is defined by
$$(\delta x)_e=\sum_{v\in V}\sigma_{v,e}x_v,$$
equal to $x_v-x_u$ when $e=\{u,v\}$ and $u$ is directed towards $v$. 

The submodule ${\rm im}(\delta)$ is the set of {\em $Q$-tensions} of
$G$. The $Q$-submodules $\ker(\partial)$ and ${\rm im}(\delta)$ are
orthogonal. 

Suppose that  $\{f_v:v\in
V\}\subset\mathbb{C}^Q$ and $\{g_e:e\in E\}\subset\mathbb{C}^Q$ are collections of functions defining
$f\in\mathbb{C}^{Q^V}$ and $g\in \mathbb{C}^{Q^E}$ by
$$f(x)=\prod_{v\in V}f_v(x_v),\hspace{1cm}g(y)=\prod_{e\in E}g_e(y_e).$$
In \cite[Theorem 12]{AJG05} it is shown that
\be\label{eqn: gen duality}q^{-\frac12|V|}\sum_{x\in Q^V}\prod_{v\in V}f_v(x_v)\prod_{e\in
  E}\ol{g_e}((\delta x)_e)=q^{-\frac12|E|}\sum_{y\in Q^E}\prod_{v\in
  V}f^F_v((\partial y)_v)\prod_{e\in E}\ol{g_e^F}(y_e),\ee
where the bar denotes complex conjugation and $F$ is the unitary
Fourier transform. In more compact notation,
$$q^{-\frac12|V|}\langle f,g\circ\delta\rangle=q^{-\frac12|E|}\langle
f^F\circ\partial,g^F\rangle.$$ 
This identity is a generalization of the Poisson summation formula (the case
$f=1_Q^{\otimes V}$) and will be useful in Section \ref{hierarchy}.

\section{Edge and vertex colouring models}\label{sec: vertex models edge}

\subsection{A general definition. Orthogonal symmetry of edge
  colouring models}

We define partition functions of vertex and edge colouring models in as great a
generality as required for this paper. 
(An example not covered by this definition is the general
$q$-state Potts model where interactions at an edge $e$ depend on $e$
as well as the colours on its endpoint vertices.\footnote{ Using the
Fortuin-Kasteleyn representation \cite{FK72} of the $q$-state Potts model, this is
 the multivariate Tutte polynomial on the hyperbola $H_q$. See e.g. \cite{Sokal}.})
\medskip

\begin{definition} \label{def: vertex edge} Let $G=(V,E)$ be a graph
  and $Q$ a set of size $q$.

A {\em partition function of a vertex $Q$-colouring model} with weight functions
$f\in\mathbb{C}^Q$ and $g\in\mathbb{C}^{Q^2}$ is a sum of the form
$$\sum_{x\in Q^V}\prod_{v\in V}f(x_v)\prod_{e\in E}g(x_v:v\in e).$$

A {\em partition function of an edge $Q$-colouring model} with weight functions
$f\in\mathbb{C}^{Q^*}$ and $g\in\mathbb{C}^Q$ is a sum of the form
$$\sum_{y\in Q^E}\prod_{v\in V}f(y_e:e\ni v)\prod_{e\in E}g(y_e).$$ 
\end{definition}

Up until Section \ref{hierarchy} we shall work only with {\em uniform} vertex
$Q$-colouring models, namely those for which $f=1_Q$, i.e. only with partition functions of the form
$$\sum_{x\in Q^V}\prod_{e\in E}g(x_v:v\in e)=(1_{\mbox{\tiny \sc
    Monochrome}}^{\otimes V},g^{\otimes E}).$$

Likewise, a {\em uniform} edge $Q$-colouring model is one for which
$g=1_Q$, i.e. with partition function of the form 
$$\sum_{y\in Q^E}\prod_{v\in V}f(y_e:e\ni v)=(f^{\otimes
  V},1_{\mbox{\tiny \sc Monochrome}}^{\otimes E}).$$ 

 The chromatic polynomial of $G$ evaluated at $q=|Q|$ (the number of proper vertex
$q$-colourings) has a uniform vertex $Q$-colouring model with $g(a,b)$
equal to $1$ if $a\neq b$ and $0$ if $a=b$. The number
of perfect matchings of $G$ has a uniform edge $\mathbb{Z}_2$-colouring model with
weight function defined by $f(a_1,\ldots, a_d)=1$ if $\#\{i:\,
a_i=1\}=1$ and $f(a_1,\ldots, a_d)=0$ otherwise. 

The weight function $f$ in the partition function of a uniform edge colouring model is not uniquely
determined. Indeed, we begin by deriving the result of Szegedy
\cite{S05}, given as Corollary \ref{orthogonal fixes mono} below, that
partition functions of uniform edge colouring models are invariant under the action of the
group of orthogonal transformations of $\mathbb{C}^Q$ (extended to
transformations of $\mathbb{C}^{Q^*}$) on the weight function $f$. \medskip

\begin{lemma}\label{UUT=I} If $U$ is an orthogonal 
transformation $\mathbb{C}^Q\rightarrow\mathbb{C}^Q$, represented as an
orthogonal matrix with
rows and columns indexed by $Q$, 
then $U\otimes U$ is an orthogonal
transformation of $\mathbb{C}^{Q^2}$ and  
$$(U\otimes U)1_{\mbox{\tiny \sc Monochrome}\cap Q^2}=1_{\mbox{\tiny \sc
    Monochrome}\cap Q^2}.$$
\end{lemma}

\begin{proof}
For an $m\times n$ matrix $A=(a_{i,j})$, ${\rm vec}(A)$ is the $mn\times 1$
vector obtained by stacking columns of $A$ one on top of the other,
with $k$th entry $a_{i,j}$ where $i=k-\lfloor\frac{k-1}{m}\rfloor m$
and $j=\lfloor\frac{k-1}{m}\rfloor+1$. 
It is easy to verify that for matrices $A,B,C$ of compatible
dimensions
$(B^T\otimes A){\rm vec}(C)={\rm vec}(ACB).$

The function $1_{\mbox{\tiny \sc Monochrome}}$ on $\mathbb{C}^{Q^2}$
when represented by a column vector indexed by $Q^2$
is equal to ${\rm vec}(I)$, where $I$ is the identity matrix with rows
and columns indexed by $Q$. 
Since $(U\otimes U){\rm vec}(I)={\rm vec}(UIU^T)={\rm vec}(I)$,
the statement of the lemma follows. 
\end{proof}

Lemma \ref{UUT=I}, along with the fact that an orthogonal
transformation preserves euclidean inner products, yields the
invariance of partition functions of uniform edge colouring models
under orthogonal transformations of the vertex weight.\medskip

\begin{corollary}\label{orthogonal fixes mono} {\rm \cite[Proposition 2.3]{S05}}
If $U$ is an orthogonal transformation of $\mathbb{C}^Q$ 
then
$$(f^{\otimes V},1_{\mbox{\tiny \sc Monochrome}}^{\otimes
  E})=((f^U)^{\otimes
  V},1_{\mbox{\tiny \sc Monochrome}}^{\otimes E}).$$
\end{corollary}

\subsection{Uniform vertex and edge colouring models: the Tutte
  polynomial and the symmetric weight enumerator of $Q$-flows}\label{Tutte}

The Tutte polynomial on the hyperbola $H_q:=\{(s,t):(s-1)(t-1)=q\}$ is
defined by
$$(s-1)^{|E|-r(E)}T\left(G;s,\frac{s-1+q}{s-1}\right)=\sum_{A\subseteq
  E}q^{|A|-r(A)}(s-1)^{|E|-|A|}.$$
Alternatively
$$(t-1)^{r(E)}T\left(G;\frac{t-1+q}{t-1},t\right)=\sum_{A\subseteq E}q^{r(E)-r(A)}(t-1)^{|A|}.$$
Given $y\in Q^E$, the Hamming weight of $y$ is defined by
$|y|:=\#\{e\in E:y_e\neq 0\}$. 

The Hamming weight enumerator of the set of $Q$-flows is a
specialization of the Tutte polynomial to the hyperbola $H_q$:
\be\label{tutte as we of flows}{\rm
  hwe}(\ker(\partial);s):=\sum_{y\in\ker(\partial)}s^{|E|-|y|}=(s-1)^{|E|-r(E)}T\left(G;s,\frac{s-1+q}{s-1}\right).\ee
Dually, the Tutte polynomial on $H_q$ is also given by
\be\label{tutte as we of tensions}{\rm hwe}({\rm im}(\delta);t):=\sum_{y\in\mbox{\rm \tiny im}(\delta)}t^{|E|-|y|}=(t-1)^{r(E)}T\left(G;\frac{t-1+q}{t-1},t\right).\ee
The latter specialization of the Tutte polynomial to $H_q$ comes in the more familiar guise of the
  monochrome (or bad colouring) polynomial (see e.g. \cite{DW93}),
\be\label{eqn: bad colouring}q^{k(G)}{\rm hwe}({\rm im}(\delta);t)=\sum_{x\in Q^V}t^{\#\{uv\in
  E:x_u=x_v\}},\ee
so called because the exponent of $t$ is equal to the number of monochromatic
edges (which are bad if proper is good) in the
vertex colouring $x$. 

Taking $X$ to be a random variable with uniform distribution on
$Q^V$, by \eqref{tutte as we of tensions} the identity
\eqref{eqn: bad colouring} gives the Tutte polynomial on $H_q$ as a
uniform vertex $Q$-colouring model,
$$(t-1)^{r(E)}T\left(G;\frac{t\!-\!1\!+\!q}{t\!-\!1},t\right)=q^{r(E)}\mathbb{E}[t^{\#\{uv\in E: X_u=X_v\}}].$$

A {\em complete weight enumerator} of a set $\mathcal{S}\subseteq Q^E$
evaluated at a function $h\in\mathbb{C}^Q$ is defined by
$${\rm cwe}(\mathcal{S}; h):=\sum_{y\in\mathcal{S}}\prod_{e\in
  E}h(y_e)=\sum_{y\in\mathcal{S}}\;\prod_{c\in  Q}h(c)^{\#\{e\in
  E:y_e=c\}}.$$
The complete weight enumerator may be regarded as a multivariate polynomial over
$\mathbb{Z}$ by letting $h$ take indeterminate values, in which case
we denote it by ${\rm cwe}(\mathcal{S};~(t_c:~c\in~ Q))$ when $h(c)=t_c$ for $c\in Q$. 
The Hamming weight enumerator ${\rm hwe}(\mathcal{S};t)$ is the specialisation
$h(0)=t$ and $h(c)=1$ for $c\neq 0$. 
If $\mathcal{S}$ is a $Q$-submodule of $Q^E$ and $\mathcal{S}^\perp$
is its orthogonal submodule, then by the Poisson summation formula
(see e.g. \cite{CCC07}), 
\be\label{eqn: Poisson}{\rm cwe}(\mathcal{S};h)=q^{-|E|/2}|\mathcal{S}|{\rm
  cwe}(\mathcal{S}^\perp;h^F),\ee
where $F$ is the unitary Fourier transform.

Since $(\ker(\partial))^\perp={\rm im}(\delta)$, by \eqref{eqn: Poisson} the complete weight enumerator of $\ker(\partial)$ evaluated at a
function $h$ corresponds to a uniform
vertex $Q$-colouring model with weight function $g(a,b)=h^F(b-a)$ depending only
on the difference $b-a$ (in other words $g(a,b)=g(a+c,b+c)$ for each $c\in Q$):
\be\label{eqn: difference only}{\rm cwe}(\ker(\partial);h)=q^{\frac12(|V|-|E|)}\sum_{x\in
  Q^V}\prod_{(u,v)\in E}h^F(x_u-x_v).\ee
Note that $h^F(x_u-x_v)=h^{FN}((\delta x)_e)$, where $N$ defined by
 $h^N(c)=h(-c)$ commutes with $F$ since $F^2=N$.
The next theorem says that if $h^N=h$ (in other words, $g(a,b)=h^F(b-a)$ has the further property that
$g(a,b)=g(b,a)$), then ${\rm cwe}(\ker(\partial);h)$ has an
edge colouring model. In particular, we will deduce an edge colouring
model for the Hamming weight enumerator (Tutte polynomial on $H_q$) as
a corollary. \medskip

\begin{theorem}\label{thm: cwe vertex edge}
Let $g\in\mathbb{C}^Q$, $F$ the unitary Fourier transform on
$\mathbb{C}^Q$ and $N$ the linear transformation defined by
$g^N(a)=g(-a)$.  
 The complete weight enumerator of the set of $Q$-flows of $G$ has the
following uniform vertex colouring model and uniform edge colouring model:
\begin{align*}{\rm cwe}(\ker(\partial);g\cdot g^N)&=q^{-|V|}\sum_{x\in Q^V}\prod_{e\in
  E}\sum_{b\in Q}\prod_{v\in e}g^F(x_v-b)\\
 & = q^{-|V|}\sum_{y\in Q^E}\prod_{v\in V}\sum_{a\in Q}\prod_{e\ni
 v}g^F(a-y_e).\end{align*}
In other words, if $X$ has a uniform distribution on $Q^V$ and $Y$ a
 uniform distribution on $Q^E$ then
$${\rm cwe}(\ker(\partial);g\cdot g^{N})=q^{|E|}\mathbb{E}[\prod_{(v,e)\in H}g^F(X_v-Y_e)].$$

Similarly, for $f\in\mathbb{C}^Q$,
$${\rm cwe}({\rm im}(\delta);f\ast f^N)=q^{|E|+r(E)}\mathbb{E}[\prod_{(v,e)\in H}f(X_v-Y_e)].$$
\end{theorem}
\begin{proof} Given a graph $G$, let
$G'=(V',E')$ to be the
$2$-stretch of $G$, defined by $V'=V\cup E$ and
$E'=\{\{v,e\}:v\in e\}$, i.e. $G'$ is obtained from $G$ by replacing
each edge of $G$ by a path of length $2$.  The edges of $G'$ are in
one-one correspondence with the half-edges of $G$. Also,
$r(E')=|E|+r(E)$.

 The set of $Q$-flows of $G'$ is in one-one
correspondence with the set of $Q$-flows of $G$ as indicated by Figure
1. %\ref{fig: flow stretch}. 
The orientation of an edge $\{v,e\}$ in $G'$
is chosen so that $v$ is directed toward $e$, as illustrated. (This
makes explicit the correspondence between direct edges $(v,e)$ of $G'$
and half-edges of $G$.) 
\setlength{\unitlength}{1cm}

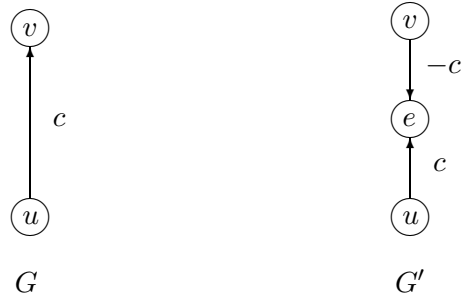
\begin{figure}[ht]\label{fig: flow stretch}
\caption{{\it Edge $e=(u,v)$ in $G$ with flow value $c$ corresponds to edges
  $(u,e)$ and $(v,e)$ in $G'$ with flow values $c$ and $-c$}}
\begin{center}
\begin{picture}(7,3.8)(0,0)
\put(0.8,0){$G$}
\put(1,1){\circle{0.5}}
\put(0.9,0.9){$u$}
\put(1,1.25){\vector(0,1){2}}
\put(1,3.5){\circle{0.5}} 
\put(0.9,3.4){$v$}

\put(1.3,2.2){$c$}

\put(5.8,0){$G'$}
\put(6,1){\circle{0.5}}
\put(5.9,0.9){$u$}
\put(6,1.25){\vector(0,1){0.8}}
\put(6,2.3){\circle{0.5}} 
\put(5.9,2.2){$e$}

\put(6,3.35){\vector(0,-1){0.8}}
\put(6,3.6){\circle{0.5}} 
\put(5.9,3.5){$v$}

\put(6.3,1.6){$c$}
\put(6.2,2.9){$-c$}

\end{picture}
\end{center}
\end{figure}
If $\mathcal{C}'$ is the set of $Q$-flows of $G'$ and $\mathcal{C}$
the set of $Q$-flows of $G$ (which is $\ker(\partial)$ in the statement of the theorem), then by the illustrated correspondence between $\mathcal{C}$ and $\mathcal{C}'$ we have
\begin{align*}{\rm cwe}(\mathcal{C};g\cdot g^N) 
&= {\rm cwe}(\mathcal{C}';g)\\
&= q^{\frac{1}{2}|E'|-r(E')}{\rm
  cwe}(\mathcal{C}'^\perp;g^F)\end{align*}
the last line by \eqref{eqn: Poisson}. With
$\frac{1}{2}|E'|-r(E')=|E|-(r(E)\!+\!|E|)=-r(E)$ and $k(G')=k(G)$, the $1$-to-$q^{k(G)}$ correspondence
between $\mathcal{C'}^\perp$ and vertex $Q$-colourings of $G'$ means
this last line can be written as the partition function of a vertex colouring model on $G'$:
\begin{align*}{\rm cwe}(\mathcal{C};g\cdot g^N)&= q^{-|V|}\sum_{z\in
  Q^{V'}}\prod_{(v,e)\in E'}g^F(z_e-z_v)\\
 &=q^{-|V|}\sum_{(x,y)\in
  Q^V\times Q^E}\prod_{(v,e)\in H}g^F(y_e-x_v)\end{align*}
using $Q^{V'}=Q^{V\cup E}\cong Q^V\times Q^E$ and $E'\cong H$ to get to the last line.

If $(X,Y)$ is uniform on $Q^V\times Q^E$, then this yields the expected
 value of the product of $g^F(X_v-Y_e)$ over half-edges $(v,e)$ as
 presented in the theorem statement. The partition functions of vertex and
 edge colouring models arise by conditioning on $Y$ and $X$ respectively.

First 
\begin{align*}\mathbb{E}[\prod_{(v,e)\in
    H}g^F(X_v-Y_e)\,\mid\, Y]&=q^{-|E|}\sum_{y\in Q^E}\prod_{(v,e)\in H}g^F(X_v-y_e)\\
 & = q^{-|E|}\prod_{e\in E}\sum_{b\in Q}\prod_{v\in e}g^F(X_v-b).\end{align*}
Hence
\begin{align*}\mathbb{E}[\prod_{(v,e)\in
    H}g^F(X_v-Y_e)] &=\mathbb{E}(\mathbb{E}[\prod_{(v,e)\in
    H}g^F(X_v-Y_e)\,\mid\, Y])\\
 & =q^{-|V|-|E|}\sum_{x\in Q^V}\prod_{e\in E}\sum_{b\in Q}\prod_{v\in e}g^F(x_v-b).\end{align*}
Second, conditioning on $X$,
\begin{align*}\mathbb{E}[\prod_{(v,e)\in
    H}g^F(X_v-Y_e)] &=\mathbb{E}(\mathbb{E}[\prod_{(v,e)\in
    H}g^F(X_v-Y_e)\,\mid\, X])\\
 & =q^{-|V|-|E|}\sum_{y\in Q^E}\prod_{v\in V}\sum_{a\in Q}\prod_{e\ni v}g^F(a-y_e).\end{align*}

For the last part of the theorem, by identity \eqref{eqn: Poisson} and $(f\ast
f^N)^F=q^{1/2}f^F\cdot f^{NF}$,
$${\rm cwe}({\rm im}(\delta);f\ast f^N)=q^{r(E)}{\rm
  cwe}(\ker(\partial); f^F\cdot f^{NF}),$$
from which the desired result follows on writing $g=f^F$ and noting that
$f^{NF}=g^{N}$ since $NF=FN$. 
\end{proof}

The Tutte
polynomial on
the hyperbola $H_q$, given by equation \eqref{tutte as we of flows}
as the Hamming weight enumerator of the set of
$Q$-flows of $G$ (for abelian group $Q$ order $q$), inherits an edge
colouring model from Theorem \ref{thm: cwe vertex edge}.\medskip

\begin{corollary}\label{cor: Tutte vertex edge}
Let $G=(V,E)$ be a graph, $Q$ a set of size $q$, and suppose $X$
takes values uniformly at random from $Q^V$ and $Y$ takes values uniformly at random from $Q^E$. 
Then
$$(s^2\!-\!1)^{|E|-r(E)}T(G;s^2,\frac{s^2\!-\!1\!+\!q}{s^2\!-\!1})=\mathbb{E}[(s\!-\!1\!+\!q)^{\#\{v\in
  e:X_v=Y_e\}}(s\!-\!1)^{\#\{v\in
  e:X_v\ne Y_e\}}].$$

In particular, the Tutte polynomial on the hyperbola $H_q$ is the
partition function of a
 uniform edge $q$-colouring model:
$$(s^2\!-\!1)^{|E|-r(E)}T(G;s^2,\frac{s^2\!-\!1\!+\!q}{s^2\!-\!1})=q^{-|E|-|V\!|}(s\!-\!1)^{2|E|}\sum_{y\in
  Q^E}\,\prod_{v\in V}\;\sum_{a\in Q}\left(\frac{s\!-\!1\!+\!q}{s\!-\!1}\right)^{\#\{e\ni v:
  y_e=a\}}.$$
\end{corollary}
\begin{proof}
Take $g=s1_0+1_{Q\setminus 0}$ in Theorem \ref{thm: cwe vertex edge},
noting that $g^N=g$ and
$g^F=q^{-\frac{1}{2}}[(s\!-\!1\!+\!q)1_0+(s\!-\!1)1_{Q\setminus 0}]$.
\end{proof}

\subsection{Real symmetric vertex colouring models as edge colouring models}\label{Vertex to Edge}

A consequence of Theorem \ref{thm: cwe vertex edge} is that if
$g\in\mathbb{C}^{Q^2}$ is symmetric, i.e. $g(a,b)=g(b,a)$, and satisfies
$g(a+c,b+c)=g(a,b)$ for each $c\in Q$, then the partition function of
a uniform vertex colouring model with weight function $g$ 
is also given by the partition function of an edge $Q$-colouring model.

If we restrict attention to real-valued vertex $Q$-colouring models
then the condition that $g\in\mathbb{R}^{Q^2}$ is symmetric is
sufficient. The key property of such a function in
this regard is
that by the spectral theorem for real symmetric matrices there is a function $h\in\mathbb{C}^{Q^2}$ such that 
\be\label{eqn: spectral decomp}g(a,b)=\sum_{c\in Q}h(a,c)h(b,c).\ee
Furthermore, if the rank of the matrix $[g(a,b)]_{a,b\in Q}$ is equal to $r$ then
there are $q-r$ values of $c$ for which $h(a,c)=0$ for all $a\in Q$. 
If the symmetric matrix $[g(a,b)]_{a,b\in Q}$ is positive semi-definite,
i.e. has non-negative eigenvalues, then the function $h$ is
real-valued. Otherwise $h$ may take purely imaginary values as well as
real values (and only purely imaginary values if the matrix
$[g(a,b)]_{a,b\in Q}$ is negative definite).
\medskip

\begin{theorem}\label{thm: szegedy} {\rm  (Szegedy \cite[Section 3.1]{S05}.)}
Suppose that $g\in\mathbb{R}^{Q^2}$ satisfies $g(a,b)=g(b,a)$ for all
$a,b\in Q$.

Then there is $h\in\mathbb{C}^{Q^2}$ such that 
 $$\sum_{x\in Q^V}\prod_{v\in V}f(x_v)\prod_{e\in E}g(x_v:v\in
 e)=\sum_{y\in  Q^E}\prod_{v\in V}\sum_{a\in Q}f(a)\prod_{e\ni v}h(a,y_e).$$

Furthermore, if the matrix $[g(a,b)]_{a,b\in Q}$ has rank $r$ then
$h(a,b)$ is identically zero for $q-r$ values of $b\in Q$ (reducing
the right-hand side to the partition function of an
edge $r$-colouring model).
\end{theorem}

\begin{proof}

Let $g(a,b)=\sum_{c\in Q}h(a,c)h(b,c)$
be a spectral decomposition of $g$.
Then
\begin{align*}\sum_{x\in Q^V}\prod_{v\in V}f(x_v)\prod_{e\in E}g(x_v:v\in e) 
& =\sum_{x\in Q^V}\prod_{v\in V}f(x_v)\prod_{e\in
  E}\sum_{c\in Q}\prod_{v\in e}h(x_v,c)\\
 & =\sum_{x\in Q^V}\prod_{v\in V}f(x_v)\sum_{y\in Q^E}\prod_{e\in E}\prod_{v\in
  e}h(x_v,y_e)\\
 & = \sum_{y\in Q^E}\sum_{x\in Q^V}\prod_{v\in V}f(x_v)\prod_{e\in E}\prod_{v\in
  e}h(x_v,y_e)\\
 & =\sum_{y\in Q^E}\prod_{v\in V}\sum_{a\in Q}f(a)\prod_{e\ni v}h(a,y_e).\end{align*}
\end{proof}

\begin{remark}\label{note: alt proof}
An alternative proof of Theorem \ref{thm: cwe vertex edge} is to
mimick the proof of Theorem \ref{thm: szegedy} applied to the case
${\rm cwe}({\rm im}(\delta);f\ast f^N)$, the latter the partition
function of a uniform vertex
$Q$-colouring model with weight function $g(a,b)$ given by
$$g(a,b)=f\ast f^N(b-a)=\sum_{c\in Q}f(c-a)f(c-b),$$
which takes the form displayed in equation \eqref{eqn: spectral decomp} with
$h(a,b)=f(b-a)$.
\end{remark}

\subsection{A family of chromatically defined graph
  polynomials}\label{hierarchy}

So far we have seen how the Tutte polynomial on $H_q$ (or monochrome
polynomial) and its generalization to
a symmetric weight enumerator of $Q$-flows (a complete weight
enumerator ${\rm cwe}(\ker(\partial);h)$ satisfying  $h^N=h$) have
uniform edge $Q$-colouring models. 

A generalization of the chromatic polynomial $P(G;q)$ in a direction
 different to that of the Tutte polynomial and complete weight enumerator
 is the symmetric function analogue of the monochrome polynomial
 \eqref{eqn: bad colouring}. This is a polynomial in commuting indeterminates $s_0,s_1,\ldots$
 defined \cite{St98} by  
\be\label{eqn: symmetric monochrome}X(G;s_0, s_1,\ldots;t)=\sum_{x\in \mathbb{N}^V}t^{\#\{uv\in
  E:x_u=x_v\}}\prod_{a\in \mathbb{N}}s_{a}^{\#\{v\in V:x_v=a\}}.\ee
This function is invariant under permutations of $s_0,s_1,\ldots$ and specializes
to Stanley's generalized symmetric chromatic function \cite{St95}
  $X(G;s_0,s_1,\ldots)$ upon setting $t=0$. 

Consider the specialization $0=s_q=s_{q+1}=\cdots$ which has the
effect of restricting the range of summation in \eqref{eqn:
  symmetric monochrome} to vertex $q$-colourings:
$$X_q(G;s_0,\ldots, s_{q-1};t):=\sum_{x\in\{0,\ldots, q-1\}^V}t^{\#\{uv\in
  E:x_u=x_v\}}\prod_{a\in \{0,1,\ldots, q-1\}}s_a^{\#\{v\in V:x_v=a\}}.$$
Note that
$$X_q(G;1,1,\ldots ,1;t)=q^{k(G)}(t-1)^{r(E)}T(G;\frac{t-1+q}{t-1},t)$$ 
is the Tutte polynomial on $H_q$, and in particular $X_q(G;1,1,\ldots,
1;0)=P(G;q)$.
More generally, for an abelian group $Q$ of order $q$, define 
\be\label{eqn: gen symmetric and cwe}X_Q(G;(s_a:a\in Q),(t_b:b\in Q))=\sum_{x\in Q^V}\prod_{a\in
  Q}s_{a}^{\#\{v\in V:x_v=a\}}\prod_{b\in Q}t_b^{\#\{e\in
  E:(\delta x)_e=b\}}.\ee
This is simultaneously a generalization of the complete weight enumerator of ${\rm
  im}(\delta)$, which has $s_a=1$ for each $a\in Q$, and of
the truncated symmetric monochrome polynomial $X_q(G;s_0,\ldots, s_{q-1};t)$,
which has $Q=\{0,\ldots, q-1\}$, $t_0=t$ and
$t_b=1$ for each $0\neq b\in Q$. Indeed, the polynomial \eqref{eqn:
  gen symmetric and cwe} may be thought of as a multivariate generating
function for pairs $(x,\delta x)\in Q^V\times Q^E$, and in Lemma
\ref{lem: general polynomial duality} below we shall see that it may
dually be regarded as a generating function for pairs $(\partial
y,y)\in Q^V\times Q^E$. 

Although the polynomial defined in \eqref{eqn: gen symmetric and cwe} is
invariant under permutations of $\{s_a: a\in Q\}$, in general it is
not invariant under permutations of $\{t_b:b\in Q\}$, may also depend on the structure of $Q$ as an
abelian group, and may further depend on the orientation of $G$. For example,
take $G$ to be the graph on two vertices joined by two parallel edges.
For $Q=\mathbb{F}_4=\{0,1,\omega,\ol{\omega}\}$, the specialization
${\rm cwe}({\rm im}(\delta);(0,t_1,t_{\omega},t_{\ol{\omega}}))$ is
equal to $t_1^2+t_{\omega}^2+t_{\ol{\omega}}^2$ for all orientations
of $G$, while for $Q=\mathbb{Z}_4=\{0,1,2,3\}$, the specialization
${\rm cwe}({\rm im}(\delta);(0,t_1,t_2,t_3))$ is equal to
$t_2^2+2t_1t_3$ for cyclic orientations of $G$ and equal to
$t_1^2+t_2^2+t_3^2$ for the remaining, acyclic orientations of $G$.
 
The relationship between the various polynomials is summarized in
Figure 2 % \ref{fig: polynomials} 
(a Hasse diagram with the ordering relation being that of
specialization).

\setlength{\unitlength}{1cm}
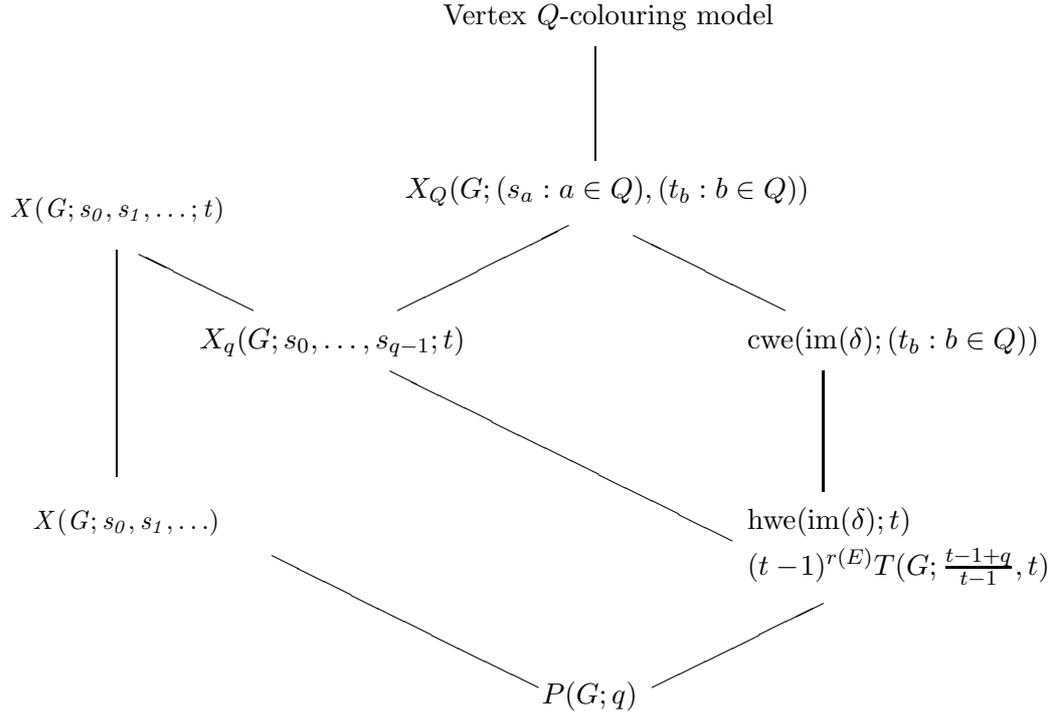
\begin{figure}[h]\label{fig: polynomials}
 \caption{{\it Some polynomials from vertex $q$-colourings of a
    graph $G$.} }
\vspace{0.5cm}
 \begin{picture}(15,10)(0,0)
\linethickness{0.5pt}

 \put(7.8,0.8){$P(G;q)$}
 \put(1.1,3.1){{\small $\mathit{X(G;s_0,s_1,\ldots)}$}}
\put(2.2,3.8){\line(0,1){3}}
\put(0.8,7.25){{\small $\mathit{X(G;s_0,s_1,\ldots ;t)}$}}
\put(7.75,1){\line(-2,1){3.5}}
\put(3.3,5.5){$X_q(G;s_0,\ldots, s_{q-1};t)$}
\put(4,6){\line(-2,1){1.5}}
\put(5.9,6){\line(2,1){2.25}}

\put(6,7.5){$X_Q(G;(s_a:a\in Q),(t_b:b\in Q))$}

\put(9.25,1){\line(2,1){2.25}}
 \put(10.5,2.5){$(t-\!1)^{r(E)}T(G;\frac{t-1+q}{t-1},t)$}
 \put(10.5,3.1){${\rm hwe}({\rm im}(\delta);t)$}
\put(10.3,2.95){\line(-2,1){4.5}}

\put(11.5,3.6){\line(0,1){1.6}}
\put(10.5,5.5){${\rm cwe}({\rm im}(\delta);(t_b:b\in Q))$}
\linethickness{0.5pt}
\put(11,6){\line(-2,1){2}}

\put(8.5,8){\line(0,1){1.5}}

\put(6.5,9.8){Vertex $Q$-colouring model}

 \end{picture}

\end{figure}

\begin{lemma}\label{lem: general polynomial duality} The partition
  function $X_Q(G;(s_a:a\in Q),(t_b:b\in Q))$ of a vertex colouring model  
whose edge weights depend only on the
  boundary $\delta$ has the following expansion over edge colourings
  with vertex weights 
  depending only on the coboundary $\partial$: 
$$X_Q(G;(s_a:a\in Q),(t_b:b\in Q))=q^{-|E|}\sum_{y\in Q^E}\prod_{a\in
  Q}\wh{s_a}^{\,\#\{v\in V:(\partial y)_v=a\}}\prod_{b\in
  Q}\wh{t_b}^{\,\#\{e\in E:y_e=b\}},$$
where 
$$\wh{s_a}=\sum_{c\in Q}\ol{\chi}(ca)s_c,\; \wh{t_b}=\sum_{c\in
  Q}\chi(cb)t_c$$
for a generating character $\chi$ of $Q$.
\end{lemma}
\begin{proof}
Apply the duality identity \eqref{eqn: gen duality} with $f_v(x_v)=s_{x_v}$ and
$\ol{g_e}(y_e)=t_{y_e}$.
\end{proof}

The polynomial $X_q(G;1,s,\ldots, s^{q-1};t)$ is called the principal
  specialization of order $q$ \cite{St99} of the symmetric monochrome
  polynomial. Loebl \cite[Theorem 3]{Loebl06} proved that the  principal
  specialization of order $2$ has an edge colouring model analogous to Van der
Waerden's eulerian subgraph expansion of the Ising model \cite{vdW41}.
  Loebl's result is the case $q=2$ of Theorem \ref{thm: van der W}
  below, Van der Waerden's expansion the case $q=2, s=1$.\medskip

  \begin{theorem} \label{thm: van der W} \hspace{0.5cm}

\begin{itemize}

\item[(i)] When $s^q\neq 1$,
$$X_q(G;1,s\ldots,s^{q-1};t)=q^{-|E|}(t-\!1)^{|E|}\sum_{y\in
    \mathbb{Z}_q^E}\left(\frac{t\!-\!1\!+\!q}{t\!-\!1}\right)^{|E|-|y|}\prod_{v\in
    V}\left(\frac{s^q-1}{se^{\frac{2\pi i(\partial y)_v}{q}}\!-\!1}\right).$$

\item[(ii)] When $s=e^{-2\pi i c/q}$ for
some $c\in\{0,\ldots, q-1\}$, %Define 
$$X_q(G;1,e^{-2\pi i c/q},\ldots, e^{-2\pi i(q-1)c/q};t)=q^{|V|-|E|}(t-\!1)^{|E|}\mathop{\sum_{y\in\mathbb{Z}_q^E}}_{\forall_{v\in V}\;(\partial y)_v=c}\left(\frac{t-1+q}{t-1}\right)^{|E|-|y|}.$$
\end{itemize}
\end{theorem}
\begin{proof} In Lemma \ref{lem: general polynomial duality} take
  $Q=\mathbb{Z}_q\cong\{0,1,\ldots, q-1\}$, $\chi(c)=e^{2\pi i c/q}$, $s_a=s^a$
  and $t_0=t, t_b=1$ for $b\neq 0$.
\end{proof}

Apart from the case $q=2$, the expansion of the principal
specialization order $q$ of Stanley's symmetric monochrome polynomial
given in Theorem \ref{thm: van der W} is not an edge colouring model
partition function, but rather a weighted sum of edge colourings with
the weight at a vertex in a given colouring depending not
only on the colours of incident edges but also on their orientation
(in order to determine the boundary of the edge colouring at the vertex).

From Theorem \ref{thm: cwe vertex edge} every symmetric weight
enumerator of $Q$-tensions has a uniform edge colouring
model. By the proof of Theorem \ref{thm: szegedy} it is easy to see
that this property extends to the non-uniform
 vertex colouring model partition function $X_Q(G;(s_a:a\in Q),(t_b:b\in
Q))$, provided of course that the edge weights
are symmetric (independent of the orientation of $G$),
i.e. $t_b=t_{-b}$.  \medskip
 
\begin{theorem}\label{thm: general chrom edge}
If $t_{-b}=t_b$ for each $b\in Q$ then
$$X_Q(G;(s_a:a\in Q),(t_b:b\in Q))=\sum_{y\in Q^E}\prod_{v\in
  V}\sum_{a\in Q}s_a\prod_{b\in Q}u_b^{\#\{e\ni v:y_e=a+b\}},$$
where $u_b$ is defined by %$u\ast u^N=t$, i.e.
 $t_b=\sum_{a\in
  Q}u_au_{a-b}$.

In particular, 
 the principal specialization order $q$ of the symmetric
monchrome function has edge colouring model given by
$$(t^2\!-\!1)^{|E|}X_q(G;1,s,\ldots,
s^{q-1};\frac{t^2\!-\!1\!+\!q}{t^2\!-\!1})$$
$$=q^{-|V\!|-\!|E|}(t\!-\!1)^{2|E|}\sum_{y\in \{0,1,\ldots, q-1\}^E}\prod_{v\in
  V}\sum_{a\in\{0,1,\ldots,
  q-1\}}s^a\left(\frac{t\!-\!1\!+\!q}{t\!-\!1}\right)^{\#\{e\ni v:y_e=a\}}.$$
\end{theorem}

\begin{proof}
Follow the proof of Theorem \ref{thm: szegedy} using the decomposition
of the weight $t_b$ as a convolution $(u\ast u^N)_b$, i.e. as a
spectral decomposition of the form needed for the proof of this
earlier theorem to go through. Cf. Remark \ref{note: alt proof}.

In order to obtain the special case, take $t_0=t^2-1+q$, $t_b=t^2-1$
for $b\neq 0$, for which $q^{1/2}u_0=t-1+q$ and $q^{1/2}u_b=t-1$ for $b\neq 0$. 
\end{proof}

\section{Non-symmetric edge colouring models for proper edge colourings}\label{Cyclic}

\subsection{Oriented (near) $2$-factorizations of $k$-regular graphs, parity and proper
  edge colourings}\label{sec: 2-factorizations}

In this section $G=(V,E)$ will be a $k$-regular graph, with line graph
$L(G)=(E,L)$ whose edge set is defined by setting $\{e,f\}\in L$ if $e$
  and $f$ are incident with a common vertex in $G$. To simplify
  exposition we shall also take
  $Q=\mathbb{Z}_q$ for some $q\geq k$, although the results of this
  section can easily be adapted to other abelian groups of order
  $q$. 

An arbitrary linear order is put on $\mathbb{Z}_q$, say the usual
integer order $0\!<\!1\!<\!2\!<\!\cdots\!<\!q\!-\!1$. A subset $K\subseteq \mathbb{Z}_q$
inherits a linear order from $\mathbb{Z}_q$ that enables us to assign a parity to a permutation of $K$, {\em
  even} (sign $+1$) if alternating, and {\em odd} (sign $-1$)
otherwise. It will be helpful to make a more general
definition.\medskip

\begin{definition}\label{def: sgn}
Given an injection $\beta:X\rightarrow Y$ between two
linearly ordered sets $(X,<)$ and $(Y,<)$, the sign of $\beta$ is
defined by 
$${\rm sgn}(\beta)=(-1)^{\#\{\ell,m\in X:\; \ell<m,\;
  \beta(\ell)>\beta(m)\}}.$$
If $\beta$ is not injective set ${\rm sgn}(\beta)=0$.
\end{definition} 

We shall also suppose that a linear order $<$ has been put on the half-edge set $H(v)=\{(v,e):\;
e\ni v\}$ for each $v\in V$. A half-edge colouring $z\in \mathbb{Z}_q^H$ with the
property that the restricted map $z_v:H(v)\rightarrow\mathbb{Z}_q$ is injective
for each $v\in V$ has sign defined by
$${\rm sgn}(z)=\prod_{v\in V}{\rm sgn}(z_v)=\prod_{v\in
  V}(-1)^{\#\{(v,e)<(v,f):\; z_{(v,e)}>z_{(v,f)}\}}.$$ A similar definition of
sign holds for an edge colouring $y\in \mathbb{Z}_q^E$, namely ${\rm sgn}(y):={\rm sgn}(z)$ where
$z\in \mathbb{Z}_q^H$ is defined by $z_{(v,e)}=y_e$ for each $(v,e)\in H$.

Signs of proper edge $k$-colourings of $G$ have been considered for
$k=3$ by many authors, the earliest (implicit) example being Vigneron
\cite{Vign46}. Penrose in
\cite{P71} first explicitly stated that the sign of a proper edge
$3$-colouring of a plane cubic graph is constant and independently Scheim \cite{S74}
proved it. Other accounts of
this property of plane cubic graphs can be found in \cite{FJ89}, \cite{K90}, \cite{MA97}. More recently, Ellingham and Goddyn
\cite{EG05} proved that a similar result held for $k$-regular plane
graphs for $k\geq 3$, and
Norine and Thomas \cite{NT07}, \cite{RT07} extended this to the larger class
of $k$-regular graphs that admit Pfaffian labellings (a generalization
of Pfaffian orientations \textemdash see below).

We use the fact that $k$-regular graphs
admitting Pfaffian labellings have proper edge $k$-colourings of
constant sign in order to construct a
function $f\in\mathbb{C}^{\mathbb{Z}_q^k}$ for which $(f^{\otimes V},1_{\mbox{\sc \tiny Zero-sum}}^{\otimes
  E})$ is up to sign equal to $P(L(G);k)$.
The following lemma then
converts this inner product into the desired edge $q$-colouring model. \medskip

\begin{lemma}\label{lem: inner zero to mono by F} 
For graph $G=(V,E)$ and $f\in\mathbb{C}^{Q^*}$,
$$(f^{\otimes V},1_{\mbox{\sc \tiny Zero-sum}}^{\otimes
  E})=((f^F)^{\otimes V},1_{\mbox{\sc \tiny Monochrome}}^{\otimes
  E}).$$
\end{lemma}
\begin{proof}
This follows since the Fourier transform $F$ is unitary, and
$(\mbox{\sc Zero-sum}\cap Q^2)^\perp=\mbox{\sc Monochrome}\cap Q^2$.
\end{proof}

A proper edge $k$-colouring of a $k$-regular graph is a partition of
the edge set into $k$ edge-disjoint perfect matchings.
The union of any two of these perfect matchings is a $2$-factor of $G$
with the property that all its components have even size.

Let $P\subseteq \mathbb{Z}_q$  have the property that $|P\cup(-P)|=k$
and that $P\cap (-P)\subseteq\{0\}$ if $q$ is odd or
$P\cap(-P)\subseteq\{0,\frac{q}{2}\}$ if $q$ is even. Thus either
$0\in P$ or $\frac{q}{2}\in P$ if $k$ is odd while either $\{0,\frac{q}{2}\}\subseteq P$
or $P\subseteq\mathbb{Z}_q\setminus\{0,\frac{q}{2}\}$ if $k$ is
even. For example, when $q=k$ we can take $P=\{0,1,\ldots,
\frac{k-1}{2}\}$ for $k$ odd and $P=\{0,1,\ldots, \frac{k}{2}\}$ for
$k$ even. 
Set $K=P\cup(-P)$, a subset of $\mathbb{Z}_q$ of size $k$. 

An {\em ordered (near)
  $2$-factorization} of $G=(V,E)$ is an ordered partition $\mathcal{F}=(F_a~:~a\in~P)$ of $E$, where each $F_a$ is a $2$-factor of $G$ if $a\neq -a$ and
  a $1$-factor of $G$ if $a=-a$. (So if $k$ is odd then there is one
  $1$-factor labelled either $F_0$
  or $F_{q/2}$, and if $k$ is even there are either no $1$-factors or $F_0$ and $F_{q/2}$ are both
  $1$-factors.) 
An {\em orientation} of $\mathcal{F}$ is an orientation of $G$
  with the property that each $2$-factor $F_a$ in $\mathcal{F}$ is a union of directed
  circuits; the $1$-factor(s) $F_0$ (and $F_{q/2}$) can be arbitrarily oriented.
In an ordered {\em bipartite} (near) $2$-factorization $(F_a:a\in P)$
  each $2$-factor $F_a$, $a\neq -a$, is a union of even length
  circuits. In other words, the $2$-factor $F_a$ is a union of two
  edge-disjoint $1$-factors.

Recall that for each $v\in V$ there is a linear order $<$ on the set
$H(v)=\{(v,e):e\ni v\}$, the $k$ half-edges incident with the vertex
$v$.
If $G$ is embedded in an orientable surface, this order
might be taken in a clockwise sense around $v$, for example. The order up to even
permutation of half-edges incident with $v$ will be used to distinguish the
bipartite ordered (near) $2$-factorizations of $G$ from the
non-bipartite.\footnote{When $k$ is odd, each way of taking the half-edges
  $\{(v,e):e\ni v\}$ incident with a vertex $v$ in a clockwise order is related to any other
  clockwise order by an even permutation (a $k$-cycle), and to
  anticlockwise orders by an odd permutation. When $k$ is even, taking the
half-edges in a clockwise direction alternates parity according to
which half-edge starts the order. Thus, for even $k$, in order to relate the linear order up to
even permutation on
$\{(v,e):e\ni v\}$ to the embedding of $G$ an additional
topological aspect of the embedding needs to be used to determine
where to start the clockwise order (the ``$0$-consistent'' order taken
in \cite{EG05} for example).} Two linear orders on $H(v)=\{(v,e):e\ni v\}$
differ by either an even permutation or an odd permutation. Likewise,
 two sets of orders on $\{H(v):v\in V\}$ differ altogether either by an even or
 odd permutation.  By the phrase `` the order up to even
permutation of half-edges around vertices'' we shall mean that
two sets of linear orders which belong to the same one of these two parity classes are equivalent.

Define the edge colouring $y\in P^E$ by setting  $y_e=a$ if $e\in F_a$. 
When the ordered (near) $2$-factorization $\mathcal{F}$ is oriented by $\sigma$,  the map
$H(v)\rightarrow K, \; (v,e)\mapsto\sigma_{v,e}y_e$ is an injection, 
and we define 
$${\rm sgn}(\mathcal{F})=\prod_{v\in V}{\rm sgn}((v,e)\mapsto \sigma_{v,e}y_e).$$

The colours $(\sigma_{v,e}y_e:e\ni v)$ are a $k$-tuple $(b_0,b_1,\ldots, b_{k-1})$
of distinct elements of $K$. Two such $k$-tuples may differ either by
an even or an odd permutation, and we define two parity classes accordingly. Let
$$\mbox{\sc Even}(K)=\{(b_0,b_1,\ldots, b_{k-1})\in K^k:\; {\rm sgn}(\ell\mapsto b_\ell)=+1\},$$
and define $\mbox{\sc  Even}\subset\mathbb{Z}_q^k$ by
$$\mbox{\sc  Even}=\bigcup_{K\subseteq \mathbb{Z}_q, |K|=k}\mbox{{\sc Even}(K)}.$$
The sets $\mbox{\sc Odd}(K)$ and $\mbox{\sc Odd}$ are defined similarly. The set
$\mbox{\sc  Even}\cup\mbox{\rm \sc Odd}$ comprises all $k$-tuples
$(b_0,b_1,\ldots, b_{k-1})$ such that $b_0,b_1,\ldots b_{k-1}$ are
distinct elements of $\mathbb{Z}_q$. \medskip

\begin{lemma}\label{lem: OOB2F} Let $G=(V,E)$ be a $k$-regular graph
  and $q\geq k$. Suppose $K\subseteq\mathbb{Z}_q$ is of size
  $|K|=k$ and satisfies $-K=K$. Then 
$$\sum_{\mathcal{F}}{\rm sgn}(\mathcal{F})=\pm ((1_{\mbox{\tiny \sc Even$(K)$}}-1_{\mbox{\tiny \sc
    Odd$(K)$}})^{\otimes V},1_{\mbox{\tiny \sc
    Zero-sum}}^{\otimes E}),$$
where the sum is over all oriented ordered bipartite (near) $2$-factorizations
 $\mathcal{F}$ of $G$ and the sign on the right-hand side depends on
 the order up to even permutation of half-edges around vertices.  
\end{lemma} 
\begin{proof} Fix an orientation $\tau$ of $G$.
Suppose we are given the
  unique representation of
  $\mathcal{F}$ as an orientation $\sigma$ cyclic on its $2$-factors
  together with an edge colouring $y\in P^E$ such that
  $\{e:y_e=a\}=F_a$. The $k$-tuple
  $(\sigma_{v,e}y_e:e\ni v)$ is then a permutation of $K$. There is a
  unique $z\in K^E$ such that $(\tau_{v,e}z_{e}:e\ni v)$ is a
  permutation of $K$, namely $z$ is the edge colouring defined by
  $z_e=y_e$ if $\tau$ preserves the direction of $\sigma$ on the edge $e$, and $z_e=-y_e$
  otherwise. (So $z_e=\tau_{v,e}\sigma_{v,e}y_e$ for either choice of vertex $v\in e$.) Hence the oriented ordered (near)
  $2$-factorizations are counted by   
$\pm ((1_{\mbox{\tiny \sc Even$(K)$}}+1_{\mbox{\tiny \sc
    Odd$(K)$}})^{\otimes V},1_{\mbox{\tiny \sc
    Zero-sum}}^{\otimes E}).$

Given $\mathcal{F}=(F_a:a\in P)$ and $y\in P^E$ defined as above by
  $y_e=a$ when $e\in F_a$, switching the orientation $\sigma$ on a
  component circuit $C$ of a $2$-factor $F_a$ corresponds
  to switching the sign of $y$ on the
  circuit $C$. Under the fixed orientation $\tau$ of $G$, the
  corresponding operation is to switch the sign of $z\in K^E$ on the
  circuit $C$, where $z$ is defined
  as above (negating $y$ when $\tau$ reverses $\sigma$).

 For a vertex $v$ belonging to $C$, the effect on the $k$-tuple
 $(\tau_{v,e}z_e:e\ni v)$ is to transpose $a$ and $-a$, thus switching
 the sign of the injection $((v,e):e\ni v)\rightarrow K,\; (v,e)\mapsto\tau_{v,e}z_e$. 

Hence, if a component
  circuit of $F_a$ has an
  odd number of vertices the factorization $\mathcal{F}$ to which it belongs will not be counted in the sum
$((1_{\mbox{\tiny \sc Even$(K)$}}-1_{\mbox{\tiny \sc
    Odd$(K)$}})^{\otimes V},1_{\mbox{\tiny \sc
    Zero-sum}}^{\otimes E})$,
and the lemma is proved. 
\end{proof}

\subsection{The Fourier transform of the parity function}\label{sec:
  fourier parity}
Our goal now is to calculate the Fourier transform of the function $1_{\mbox{\tiny \sc Even$(K)$}}-1_{\mbox{\tiny \sc
    Odd$(K)$}}$ for $K\subseteq \mathbb{Z}_q$ of size $k$ satisfying $-K=K$. First
    we shall do so for $q=k$ and $K=\mathbb{Z}_k$, thereby establishing
    that the sum in Lemma \ref{lem: OOB2F} is also equal to a
    sum over proper edge $k$-colourings of $G$ (Corollaries  \ref{cor: zero sum to mono 2} and \ref{cor: OOB2F to EC} below). \medskip

\begin{lemma}\label{lem: det F}
The matrix $q^{1/2}F=[e^{2\pi i\ell m/q}]_{0\leq \ell, m\leq q-1}$ has
determinant 
$$\det[e^{2\pi i\ell m/q}]_{0\leq \ell, m\leq
  q-1}=i^{(q-1)(3q-2)/2}q^{q/2}.$$
\end{lemma}
\begin{proof}
 Since $F^2=N$ (recall that $N=[1_\ell(-m)]_{0\leq \ell,m\leq q-1}$), 
$$\det((q^{1/2}F)^2)=q^q\det(N)=q^q(-1)^{q-1}$$
so that $\det[e^{2\pi i\ell m/q}]_{0\leq \ell, m\leq q-1}=\pm
i^{q-1}q^{q/2}$. To determine the sign, note that
$$e^{2\pi i m/q}-e^{2\pi i\ell/q}=2i(e^{\pi
  i/q})^{\ell+m}\sin\frac{(m-\ell)\pi}{q},$$
and
$$\sum_{0\leq \ell<m\leq q-1}\ell+m=\frac{q(q-1)^2}{2}.$$
Using the well known expansion of Vandermonde determinants, we obtain
\begin{align*}\det[e^{2\pi i\ell m/q}]_{0\leq \ell, m\leq q-1}
&=\prod_{0\leq \ell<m\leq q-1}(e^{2\pi i m/q}-e^{2\pi i\ell/q})\\
& = (e^{\pi i/q})^{q(q-1)^2/2}i^{q(q-1)/2}\prod_{0\leq \ell<m\leq
  q-1}2\sin\frac{(m-\ell)\pi}{q},\end{align*}
and $\sin\frac{\pi(m-\ell)}{q}>0$ for $0\leq \ell<m\leq q-1$. 
This yields the required sign and the expression for the determinant given in the statement of the
lemma now follows.
\end{proof}

\begin{corollary}\label{cor: zero sum to mono 2} For a $k$-regular
  graph $G=(V,E)$,
$$((1_{\mbox{\tiny \sc Even$(\mathbb{Z}_k)$}}-1_{\mbox{\tiny \sc
    Odd$(\mathbb{Z}_k)$}})^{\otimes V},1_{\mbox{\tiny \sc Zero-sum}}^{\otimes E})=\pm((1_{\mbox{\tiny \sc Even$(\mathbb{Z}_k)$}}-1_{\mbox{\tiny \sc
    Odd$(\mathbb{Z}_k)$}})^{\otimes V},1_{\mbox{\tiny \sc
    Monochrome}}^{\otimes E}),$$
where the sign is given by
$$\begin{cases}(-1)^{\frac{k-1}{2}|E|} &  \mbox{\rm $k$ odd,}\\
 (-1)^{\frac{k}{2}|E|+\frac{|V|-|E|}{2}} & \mbox{\rm $k$ even.}\end{cases}$$
\end{corollary}
\begin{proof} With a view to applying Lemma \ref{lem: inner zero to mono by F},
  we calculate $(1_{\mbox{\tiny \sc Even$(\mathbb{Z}_k)$}}-1_{\mbox{\tiny \sc
    Odd$(\mathbb{Z}_k)$}})^F$.
  
For $(b_0,\ldots, b_{k-1})\in \mathbb{Z}_k^k$,
the matrix $B:=[k^{-1/2}e^{\frac{2\pi i\ell b_\ell}{k}}]_{0\leq
  \ell\leq k-1}$ is zero if $b_\ell=b_m$ for some $0\leq \ell<m\leq
k-1$. Otherwise $(b_0,b_1\ldots, b_{k-1})$ is a permutation $\beta$ of
  $(0,1,\ldots, k-1)$ and the matrix $B$ is equal to the Fourier matrix $F=k^{-1/2}[e^{2\pi i\ell
  m/k}]_{0\leq\ell, m\leq k-1}$ with its columns permuted
  by $\beta$, so that by Lemma \ref{lem: det F}
$$\det[k^{-1/2}e^{2\pi i\ell b_\ell/q}]_{0\leq \ell\leq q-1}={\rm sgn}(\beta)i^{(3q-2)(q-1)/2}.$$
Hence
\begin{align*}(1_{\mbox{\tiny \sc Even$(\mathbb{Z}_k)$}}-1_{\mbox{\tiny \sc
    Odd$(\mathbb{Z}_k)$}})^F(b_0,\ldots, b_{k-1})
&=k^{-k/2}\sum_{\rho\in{\rm
    Sym}\{0,1,\ldots, k-1\}}{\rm sgn}(\rho)\prod_{0\leq \ell\leq
    k-1}e^{2\pi i\rho(\ell)b_\ell/q}\\
&=i^{(k-1)(3k-2)/2}(1_{\mbox{\tiny \sc Even$(\mathbb{Z}_k)$}}-1_{\mbox{\tiny \sc Odd$(\mathbb{Z}_k)$}}).\end{align*}
With $(k-1)(3k-2)/2\equiv (k^2-k+1+(-1)^k)/2\,(\mbox{\rm mod}\, 4)$,
and $k|V|=2|E|$ since $G$ $k$-regular, when taking the product over
all vertices this makes the exponent of $i$ equal to
$(k-1)|E|+|V|$ when $k$ is even and $(k-1)|E|$ when $k$ is odd.
Note that although $(k-1)|E|+|V|$ can be odd when $|V|$ is odd, in
this case the graph
$G=(V,E)$ does not have a proper edge $k$-colouring since it cannot have a
$1$-factor, i.e. the inner products in the statement of the lemma are
both zero.  
\end{proof}

From Corollary \ref{cor: zero sum to mono 2} and Lemma \ref{lem:
  OOB2F} (with $q=k, K=\mathbb{Z}_k$)
we deduce the existence of a sign-preserving bijection, exhibited by
  Alon and Tarsi in \cite{AT92}, \cite{NA93}, between ordered oriented
  bipartite (near) $2$-factorizations of $G$ and proper edge
  $k$-colourings of $G$. \medskip

\begin{corollary}\label{cor: OOB2F to EC}  Let $G=(V,E)$ be a $k$-regular graph
  and $q\geq k$. Suppose $K\subseteq\mathbb{Z}_q$ is of size
  $|K|=k$ and satisfies $-K=K$. Then 
\be\label{eqn: OOB2f to EC}\sum_{\mathcal{F}}{\rm sgn}(\mathcal{F})=\pm ((1_{\mbox{\tiny \sc Even$(K)$}}-1_{\mbox{\tiny \sc
    Odd$(K)$}})^{\otimes V},1_{\mbox{\tiny \sc
    Monochrome}}^{\otimes E}).\ee
where the sum is over all oriented ordered bipartite (near) $2$-factorizations
 $\mathcal{F}$ of $G$ and the sign on the right-hand side depends on
 the order up to even permutation of half-edges around vertices.
\end{corollary}

Corollary \ref{cor: OOB2F to EC} enables the derivation  of an edge
 $q$-colouring model for $P(L(G);k)$ whenever $G$ is a $k$-regular
 graph with a so-called a {\em Pfaffian labelling} \cite{NT07}. For a definition and discussion of Pfaffian orientations and labellings of a graph see \cite{RT07}.  All planar graphs have a Pfaffian orientation. 

A graph admits a Pfaffian $\mathbb{Z}_2$-labelling if and only if it has a
Pfaffian orientation. The complete bipartite graph $K_{3,3}$ and the Petersen graph are not Pfaffian, but the Petersen
graph admits a Pfaffian $\mathbb{Z}_4$-labelling.

In this context all that needs to be known
 about $k$-regular graphs with Pfaffian labellings is that their proper edge $k$-colourings all have the same
 sign \cite[Theorem 8.3]{RT07}, i.e. the right-hand side of equation \eqref{eqn: OOB2f to EC} is up to sign equal to $P(L(G);k)$.

From Lemma \ref{lem: OOB2F} and Corollary \ref{cor: OOB2F to
  EC} we deduce the following.\medskip

\begin{corollary} \label{lem: pfaff} 
Let $G=(V,E)$ be a $k$-regular graph that admits a
  Pfaffian labelling. 
Suppose $K\subseteq\mathbb{Z}_q$ has size $|K|=k$
  and that $-K=K$. 
Then  
\be\label{eqn: same parity edge col}((1_{\mbox{\tiny \sc
  Even$(K)$}}-1_{\mbox{\tiny \sc Odd$(K)$}})^{\otimes V},1_{\mbox{\tiny \sc Zero-sum}}^{\otimes
  E})=\pm P(L(G);k),\ee
where the sign on the right-hand side depends on
 the order up to even permutation of half-edges around vertices.  
\end{corollary}

As an illustration of Corollary \ref{lem: pfaff}, if $G=(V,E)$ is a
plane $3$-regular graph and the order of half-edges $\{(v,e):e\ni v\}$ around a
vertex $v$ up to cyclic permutation is clockwise in the plane then 
\cite{S74}, \cite{P71}, \cite{Vign46} the
sign of a proper edge $k$-colouring is always $(-1)^{|E|}$, so that
the sign in equation \eqref{eqn: same parity edge col} is $+1$.  
In \cite{EG05}, it is shown more generally that when $k\geq 3$ is odd,
a $k$-regular
plane graph with a clockwise cyclic order of half-edges around vertices always has a
positive sign for its oriented ordered near $2$-factorizations,
which by Corollary \ref{cor: zero sum to mono 2} implies that the sign
of a proper edge $k$-colouring is always $(-1)^{\frac{k-1}{2}|E|}$.
(A similar statement \cite[Lemma 3.6]{EG05} can be given for even $k$,
only now, with $k$-cycles odd permutations, a clockwise order of
half-edges needs to be consistently rooted by reference to the
embedding of $G$ in order to ensure that all oriented ordered
$2$-factorizations of $G$ have positive sign.) 

In the next lemma and its special case, Corollary \ref{cor: k+1}, we prepare for the main theorem of this section, Theorem \ref{thm:
  PLGk even odd}, which gives an
edge $q$-colouring model for $P(L(G);k)$ for $q\geq k$ when $G$ has a
  Pfaffian labelling.\medskip

\begin{lemma}\label{lem: fourier even minus odd} 
 Let $k$ be odd, $q>k$, and set $K=\{0,\pm 1,\ldots,
  \pm\frac{k-1}{2}\}\subset\mathbb{Z}_q$.
Then, for $(b_0,\ldots, b_{k-1})\in \mathbb{Z}_q^k$,
$$(1_{\mbox{\tiny \sc Even}(K)}-1_{\mbox{\tiny \sc
    Odd}(K)})^F(b_0,\ldots, b_{k-1})=q^{-k/2}i^{k(k-1)/2}\prod_{0\leq \ell<m\leq k-1}2\sin\frac{\pi
    (b_m-b_\ell)}{q}.$$
Similarly, when $k$ is even and $K=\{\pm 1,\ldots,\pm
\frac{k}{2}\}\subset \mathbb{Z}_q$, then
 for $(b_0,\ldots, b_{k-1})\in \mathbb{Z}_q^k$,
$$(1_{\mbox{\tiny \sc Even}(K)}-1_{\mbox{\tiny \sc
    Odd}(K)})^F(b_0,\ldots, b_{k-1})$$
$$=q^{-k/2}i^{k(k-1)/2}\prod_{0\leq \ell<m\leq k-1}2\sin\frac{\pi
    (b_m-b_\ell)}{q}\mathop{\sum_{S\subset\{0,1,\ldots,
      k-1\}}}_{|S|=k/2}\cos\frac{\pi(\sum_{m\in
      S}b_m-\sum_{\ell\not\in S}b_\ell)}{q}.$$

\end{lemma}
{\sc Note.} To avoid ambiguities in sign in the equations given in this lemma
assume that $b_0,b_1,\ldots, b_{k-1}$ take integer values from some
fixed system of residue classes modulo $q$, say $\{0,1,\ldots, q-1\}$.

\begin{proof} Set $f=1_{\mbox{\tiny \sc Even$(K)$}}-1_{\mbox{\tiny \sc
    Odd$(K)$}}$. 

First, take $k$ odd and $K=\{0,\pm 1,\ldots, \pm\frac{k-1}{2}\}$. 
For $(b_0,b_1,\ldots,b_{k-1})\in \mathbb{Z}_q^k$,
\begin{align*}
f^F(b_0,\ldots, b_{k-1}) & = q^{-k/2}\sum_{\rho\in \mbox{\tiny \rm
    Sym}\{-\frac{k-1}{2},\ldots, 0, \ldots, \frac{k-1}{2}\}}{\rm sgn}(\rho)\prod_{-\frac{k-1}{2}\leq\ell\leq\frac{k-1}{2}}e^{2\pi i \rho(\ell)b_{\ell+(k-1)/2}/q}\\
& = q^{-k/2}\det[e^{2\pi i m b_{\ell+(k-1)/2}/q}]_{-\frac{k-1}{2}\leq\ell,m\leq\frac{k-1}{2}}\\
& = q^{-k/2}e^{-\pi i (k-1)(b_0+\cdots b_{k-1})/q}\prod_{0\leq \ell<m\leq k-1}(e^{2\pi i b_m/q}-e^{2\pi i b_\ell/q})
\end{align*}
where the range of $\ell, m$ has been translated to $0,1,\ldots, k-1$
in the last line % by multiplying column $j$ of the matrix by
                  % $e^{2\pi i\frac{k-1}{2}b_{j}/q}$
and the resulting Vandermonde
determinant evaluated using $$\det[x_\ell^m]_{0\leq\ell, m\leq
  k-1}=\prod_{0\leq\ell< m\leq k-1}(x_m-x_\ell).$$
With
$$e^{2\pi i b_m/q}-e^{2\pi i b_\ell/q}= 2ie^{\pi(b_\ell+b_m)/q}\sin\frac{\pi (b_m-b_\ell)}{q},$$
and 
$$\sum_{0\leq\ell <m\leq k-1}b_\ell+b_m=(k-1)\sum_{0\leq \ell\leq k-1}b_\ell$$
this yields
$$f^F(b_0,b_1,\ldots,b_{k-1})=q^{-k/2}i^{k(k-1)/2}\prod_{0\leq
  \ell<m\leq k-1}2\sin\frac{\pi (b_m-b_\ell)}{q}.$$

Second, take $k$ even, $K=\{\pm 1,\pm 2,\ldots, \pm\frac{k}{2}\}\subset\mathbb{Z}_q$.
The matrix
$$[e^{2\pi i
    mb_{\ell}/q}]_{\mbox{\tiny \shortstack{$\ell\in\{0,1,\ldots,
        k-1\}$\\$m\in\{\pm 1,\ldots, \pm k/2\}$}}}$$
takes the form 
$$[x_\ell^{m-k/2}]_{\mbox{\tiny \shortstack{$\ell\in\{0,1,\ldots, k-1\}$\\$m\in\{0,1,\ldots, k\}\setminus\{k/2\}$}}}$$
where $x_\ell=e^{2\pi i b_\ell/q}$.
The determinant of the matrix $[x_\ell^m]$ (with $m\in\{0,1,\ldots,
k\}\setminus\{k/2\}$) has total degree $0+1+\cdots +k\, -k/2=k^2/2$
and is divisible by $x_m-x_\ell$ for each $0\leq
\ell<m\leq k-1$; the remaining factor is a
  homogeneous symmetric polynomial of total degree $k/2$ in
  $x_0,x_1,\ldots, x_{k-1}$, and we find that
  $$\det[x_\ell^{m}]_{\mbox{\tiny \shortstack{$\ell\in\{0,1,\ldots,
        k-1\}$\\$m\in\{\pm 1,\ldots,
        \pm k/2\}$}}}=(x_0x_1\cdots
  x_{k-1})^{-k/2}\prod_{0\leq \ell<m\leq
    q-1}(x_m-x_\ell)\mathop{\sum_{S\subset\{0,1,\ldots,
      k-1\}}}_{|S|=k/2}\prod_{\ell\in S}x_\ell.$$
Here $(x_0x_1\cdots
x_{k-1})^{-k/2}=e^{-\pi i k(b_0+\cdots +b_{k-1})/q}$.
As was seen above, 
$$\prod_{0\leq \ell<m\leq k-1}(e^{2\pi i b_m/q}-e^{2\pi i
  b_\ell/q})=e^{\pi i(k-1)(b_0+\cdots +b_{k-1})/q}\prod_{0\leq
  \ell<m\leq k-1}2\sin\frac{\pi(b_m-b_\ell)}{q},$$
so that the required determinant is given by
$$e^{-\pi i (b_0+\cdots +b_{k-1})/q}\prod_{0\leq
  \ell<m\leq k-1}2\sin\frac{\pi(b_m-b_\ell)}{q}\mathop{\sum_{S\subset\{0,1,\ldots,
      k-1\}}}_{|S|=k/2}e^{2\pi i(\sum_{m\in S}b_m)/q}$$
and the equation given in the lemma statement follows as a result.
\end{proof}

\begin{corollary}\label{cor: k+1}
(i) When $k$ is odd and $K=\mathbb{Z}_{k+1}\setminus\{\frac{k+1}{2}\}$, 
$$(1_{\mbox{\tiny \sc Even(K)}}-1_{\mbox{\tiny \sc
    Odd(K)}})^F=(k+1)^{-1/2}i^{k(k-1)/2}(1_{\mbox{\sc \tiny Even}}-1_{\mbox{\tiny \sc Odd}}).$$
(ii) When $k$ is even and $K=\mathbb{Z}_{k+1}\setminus \{0\}$,
$$(1_{\mbox{\tiny \sc Even(K)}}-1_{\mbox{\tiny \sc
    Odd(K)}})^F=(k+1)^{-1/2}i^{k(k+1)/2}(1_{\mbox{\sc \tiny Even}}-1_{\mbox{\tiny \sc Odd}}).$$
\end{corollary}

\begin{proof}
Identify $\mathbb{Z}_{k+1}$ with $\{0,1,\ldots, k\}$,  with order
$0<1<\cdots <k$, and suppose that $\beta:\ell\mapsto b_\ell$ is an
bijection from $\{0,1,\ldots, k-1\}$ to $\{0,1,\ldots, k\}\setminus\{b\}$.

(i) First consider $k$ odd. By Lemma \ref{lem: fourier even minus odd}, when $K=\mathbb{Z}_{k+1}\setminus\{\frac{k+1}{2}\}$,
$$(1_{\mbox{\tiny \sc Even}(K)}-1_{\mbox{\tiny \sc
      Odd}(K)})^F(b_0,\ldots, b_{k-1})=(k+1)^{-k/2}\prod_{0\leq \ell<m\leq k-1}2i\sin\frac{\pi
    (b_m-b_\ell)}{k+1}.$$
We have
$$ \prod_{0\leq \ell<m\leq k-1}2i\sin\frac{\pi
    (b_m-b_\ell)}{k+1}={\rm sgn}(\beta)\prod_{0\leq \ell<m\leq k}2i\sin\frac{\pi
    (m-\ell)}{k+1}/(-1)^b\mathop{\prod_{0\leq\ell\leq k}}_{\ell\neq b}2i\sin\frac{\pi(\ell-b)}{k+1}.$$
 The product in the numerator is given by Lemma \ref{lem: det F} as
$$\prod_{0\leq \ell<m\leq k}2i\sin\frac{\pi
    (m-\ell)}{k+1}=i^{k(k+1)/2}(k+1)^{\frac{k+1}{2}},$$
and the product in the denominator by
\begin{align*}\mathop{\prod_{0\leq\ell\leq k}}_{\ell\neq b}2i\sin\frac{(\ell-b)\pi}{k+1}
&=(-1)^{b}\prod_{1\leq m\leq k}e^{-\pi im/(k+1)}(e^{2\pi im/(k+1)}-1)\\
& =(-1)^{b}e^{-\pi i k(k+1)/2(k+1)}(-1)^k(k+1)\\
& = (-1)^bi^{k}(k+1),\end{align*}
using the identity $\prod_{1\leq\ell\leq q-1}(1-e^{2\pi i \ell/q})=q$
(the sum of the coefficients in the polynomial $1+t+\cdots  +t^{q-1}$
with roots $e^{2\pi i\ell/q}$).
Hence
$$(1_{\mbox{\tiny \sc Even}(K)}-1_{\mbox{\tiny \sc
      Odd}(K)})^F(b_0,\ldots, b_{k-1})={\rm
    sgn}(\beta)i^{k(k-1)/2}(k+1)^{-1/2}.$$

(ii) Second consider $k$ even, $K=\mathbb{Z}_{k+1}\setminus \{0\}$. By
the calculation for odd $k$ we see that
$$(1_{\mbox{\tiny \sc Even}(K)}-1_{\mbox{\tiny \sc
      Odd}(K)})^F(b_0,b_1,\ldots, b_{k-1})$$
$$= {\rm
    sgn}(\beta)i^{k(k-1)/2}(k+1)^{-1/2}\mathop{\sum_{S\subset\{0,1,\ldots,
      k-1\}}}_{|S|=k/2}e^{\pi i(\sum_{m\in S}b_m-\sum_{\ell\not\in S}b_\ell)/(k+1)}.$$
Since 
$$\sum_{m\in S}b_m-\sum_{\ell\not\in S}b_\ell=\sum_{m\in
  S}(b_m-b)-\sum_{\ell\not\in S}(b_\ell-b)$$
it suffices to consider the case where $b=0$, i.e. $\{b_0,b_1,\ldots,
b_{k-1}\}=\{1,2,\ldots, k\}$, in which case $\sum_{m\in
  S}b_m-\sum_{\ell\not\in S}b_\ell\equiv 2\sum_{m\in S}b_m\;(\mbox{\rm
  mod}\, k\!+\!1)$. Thus we wish to evaluate
$$\mathop{\sum_{A\subset\{1,2,\ldots, k\}}}_{|A|=k/2}e^{2\pi
  i(\sum_{a\in A}a)/(k+1)}.$$
This sum is equal to the evaluation at $s=e^{2\pi i/(k+1)}$ of the
coefficient of $t^{k/2}$ in the generating functions for subset sums:
$$[t^{k/2}]\prod_{a\in\{1,2,\ldots, k\}}(1+s^at)\;(\mbox{\rm mod}\,
s^{k+1}-1)$$
$$=\sum_{c\in\{0,1,\ldots, k\}}\#\{A\subset\{1,2,\ldots,
k\}:\; |A|=k/2,\, \sum_{a\in A}a=c\,(\mbox{\rm mod}\,
k\!+\!1)\}\,s^c.$$
With
\begin{align*}\prod_{a\in\{1,2,\ldots, k\}}(1+e^{2\pi i
  a/(k+1)}t) 
&=\prod_{a\in\{1,2,\ldots, k\}}(-t-e^{-2\pi i a/(k+1)})\\
& =\frac{t^{k+1}+1}{t+1},\end{align*}
picking out the coefficient of $t^{k/2}$ yields
$$\mathop{\sum_{A\subset\{1,2,\ldots, k\}}}_{|A|=k/2}e^{2\pi
  i(\sum_{a\in A}a)/(k+1)}=(-1)^{k/2}.$$
Hence
$$(1_{\mbox{\tiny \sc Even}(K)}-1_{\mbox{\tiny \sc
      Odd}(K)})^F(b_0,b_1,\ldots, b_{k-1})={\rm
    sgn}(\beta)i^{k(k-1)/2}(k+1)^{-1/2}(-1)^{k/2}.$$
\end{proof}

\subsection{Proper $k$-colourings by proper $q$-colourings ($q\geq
  k$)}\label{sec: last}
Finally, we can now write down the edge $q$-colouring model for proper
edge $k$-colourings of $k$-regular graphs admitting Pfaffian labellings.\medskip

\begin{theorem}\label{thm: PLGk even odd} Let $k$ be odd and $G=(V,E)$ a $k$-regular graph that admits a
  Pfaffian labelling. For each $v\in V$, suppose the half-edges 
  $\{(v,e):e\ni v\}$ around $v$ are linearly $<$-ordered. Then for $q\geq k$ the number of proper $k$-colourings of
  $G$ is the partition function of a uniform edge $q$-colouring model, given by
$$\pm P(L(G);k)=q^{-|E|}\sum_{y\in\{0,1,\ldots, q-1\}^E}\prod_{v\in V}\mathop{\prod_{e,f\ni
    v}}_{(v,e)<(v,f)}2\sin\frac{\pi (y_f-y_e)}{q},$$
where the sign depends on the order up to even permutation of
half-edges around vertices.
\end{theorem}
\begin{proof} 
Piece together Corollary \ref{lem: pfaff} and Lemma \ref{lem: fourier
  even minus odd}. 
\end{proof}

Theorem \ref{thm: PLGk even odd} has an alternative expression in
  terms of vertex $q$-colouring models. Suppose the line graph
  $L(G)=(E,L)$ of the $k$-regular graph $G$ ($k$ odd) has an
  orientation of each edge $\{e,f\}$ in $L$ defined by directing $e$
  towards $f$ whenever $(v,e)<(v,f)$ as half-edges, writing
  $(e,f)$ for this edge orientation. Then, for all $q\geq k$,
$$\pm P(L(G);k)=q^{-|E|}\sum_{y\in \{0,1,\ldots,
  q-1\}^E}\prod_{(e,f)\in L}2\sin\frac{\pi(y_f-y_e)}{q},$$
where the choice of sign depends on the orientation of $L(G)$ (up to an
even number reversals of direction on edges).
%(-1)^{\frac{k-1}{2}|E|} 

There is of course a parallel statement to  Theorem \ref{thm: PLGk
  even odd} for even $k$, only the edge colouring model has a more
complicated vertex weight (given in Lemma \ref{lem: fourier even
  minus odd}).

When $q=k$ Theorem \ref{thm: PLGk even odd} is tautologous. When
$q=k+1$ the edge $(k+1)$-colouring model for proper edge
$k$-colourings is particularly simple, and in this case both for odd
and even values of $k$ by using the result of Corollary
\ref{cor: k+1}:\medskip
 
\begin{theorem}\label{thm: k+1}
Let $G=(V,E)$ be a $k$-regular graph that admits a Pfaffian
labelling. For each $v\in V$, suppose the half-edges 
  $\{(v,e):e\ni v\}$ around $v$ are linearly ordered.
Then
$$\pm P(L(G);k)=(k+1)^{-\frac{|V|}{2}}\sum_{y\in
  \mathbb{Z}_{k+1}^E}{\rm sgn}(y),$$
where ${\rm sgn}(y)=+1$ if there are an even number of vertices $v\in
V$ such that the injection $((v,e):e\ni v)\rightarrow 0\!<\!1\!<\cdots\!<k,\;
(v,e)\mapsto y_e$ has sign $-1$, and ${\rm sgn}(y)=-1$ it there are an
odd number.
The sign on the left-hand side depends on the order up to even
permutation of half-edges around vertices.
\end{theorem}
A statement of this theorem for plane cubic graphs was given as
Proposition \ref{prop: cubic even odd 4}.

{\small 
\bibliography{bb}
\bibliographystyle{abbrv}
}

\end{document}